\documentclass[12pt,a4paper,twoside,leqno]{article}
\usepackage{epsfig}
\usepackage{graphpap,latexsym}
\usepackage{amsthm,amsmath,amssymb}
\usepackage{color}

\definecolor{grr}{rgb}{0.1,0.5,0.1}

\def\pier{}

\def\red#1{#1}
\def\blu#1{#1}
\def\grr#1{#1}

\def\trait #1 #2 #3 {\vrule width #1pt height #2pt depth #3pt}
\def\fin{
      \trait .3 5 0
      \trait 5 .3 0
      \kern-5pt
      \trait 5 5 -4.7
      \trait 0.3 5 0
\medskip}
\newtheorem{teor}{Theorem}[section]
\newtheorem{defin}[teor]{Definition}
\newtheorem{lemm}[teor]{Lemma}
\newtheorem{osse}[teor]{Remark}
\newtheorem{prop}[teor]{Proposition}
\newtheorem{defi}[teor]{Definition}
\newtheorem{coro}[teor]{Corollary}
\newtheorem{prob}[teor]{Problem}
\newtheorem{hypo}[teor]{Hypothesis}

\newcommand{\bele}{\begin{lemm}\begin{sl}}
\newcommand{\enle}{\end{sl}\end{lemm}}
\newcommand{\bedef}{\begin{defi}\begin{sl}}
\newcommand{\eddef}{\end{sl}\end{defi}}
\newcommand{\bete}{\begin{teor}\begin{sl}}
\newcommand{\ente}{\end{sl}\end{teor}}
\newcommand{\beos}{\begin{osse}\begin{rm}}
\newcommand{\eddos}{\end{rm}\end{osse}}
\newcommand{\bepr}{\begin{prop}\begin{sl}}
\newcommand{\empr}{\end{sl}\end{prop}}
\newcommand{\bepro}{\begin{prob}\begin{rm}}
\newcommand{\empro}{\end{rm}\end{prob}}
\newcommand{\bede}{\begin{defin}\begin{sl}}
\newcommand{\edde}{\end{sl}\end{defin}}
\newcommand{\beco}{\begin{coro}\begin{sl}}
\newcommand{\enco}{\end{sl}\end{coro}}
\newcommand{\behy}{\begin{hypo}\begin{sl}}
\newcommand{\enhy}{\end{sl}\end{hypo}}

\newcommand{\thspace}{\hspace{3mm}}

\newcommand{\RR}{\mathbb{R}}
\newcommand{\NN}{\mathbb{N}}

\newcommand{\beeq}[1]{\begin{equation}\label{#1}}
\newcommand{\eddeq}{\end{equation}}

\newcommand{\beeqa}[1]{\begin{eqnarray}\label{#1}}
\newcommand{\eddeqa}{\end{eqnarray}}

\newcommand{\beal}[1]{\begin{align}\label{#1}}
\newcommand{\eddal}{\end{align}}

\newcommand{\bespl}[1]{\begin{split}\label{#1}}
\newcommand{\edspl}{\end{split}}

\newcommand{\bega}[1]{\begin{gather}\label{#1}}
\newcommand{\edga}{\end{gather}}

\newcommand{\beeqax}{\begin{eqnarray*}}
\newcommand{\eddeqax}{\end{eqnarray*}}

\def\qed{\ifmmode 
  \else \leavevmode\unskip\penalty9999 \hbox{}\nobreak\hfill
  \fi
  \quad\hbox{\hskip.5em\vrule width.4em height.6em depth.05em\hskip.1em}}
\def\endproofsym{\qed}
\renewenvironment{proof}[1][Proof]{\trivlist\item[\hskip\labelsep{\hskip0pt
    {\normalfont\scshape#1.}\hskip .321429\parindent}]\ignorespaces}
{\endproofsym\endtrivlist}
\def\endnobox{\def\endproofsym{}\end{proof}\def\endproofsym{\qed}}

\newcommand{\no}{\nonumber}

\newcommand{\beeqao}{\begin{eqnarray}\no}
\newcommand{\bealo}{\begin{align}\no}
\newcommand{\besplo}{\begin{split}\no}
\newcommand{\begao}{\begin{gather}\no}

\newcommand{\duav}[1]{\langle{#1}\rangle}

\newcommand{\duavg}[1]{\left\langle{#1}\right\rangle}

\newcommand{\dt}{\partial_t}

\newcommand{\itt}{\int_0^t}

\newcommand{\io}{\int_\Omega}

\newcommand{\for}{\ \mbox{for } \,}

\newcommand{\e}{\varepsilon}
\newcommand{\ee}{_{\varepsilon}}

\renewcommand\div{\mathop{\rm div}\nolimits}

\DeclareMathOperator{\dive}{div}

\newcommand{\LDH}{L^2(0,T;H)}
\newcommand{\LDV}{L^2(0,T;V)}

\newcommand{\LIV}{L^\infty(0,T;V)}

\newcommand{\teta}{\theta}

\newcommand{\nat}{\mathbb{N}}

\def\vp{\varphi}
\def\dvp{\partial\varphi}
\def\domvp{{\mathcal D}(\varphi)}
\def\dcvp{{\mathcal D}_C(\varphi)}
\def\dczvp{{\mathcal D}_{C_0}(\varphi)}

\def\vrh{\varrho}

\newfont{\ctv}{msam10}
\newcommand{\bbox}{\mbox{\ctv \symbol{4}}}
\def\QED{{$\hfill\bbox$}}

\def\tk{2^j}
\def\tko{2^{j-1}}
\def\pk{\Phi_j}
\def\tkm{2^{-j}}

\def\vr{\varrho}

\newcommand{\iop}{\int_{\Omega'}}
\newcommand{\iopd}{\int_{\partial\Omega'}}

\def\fine{\hfill\kern4pt \vrule height4pt depth0pt width4pt }

\def\real{\mathbb{R}}
\def\for{\mbox{ \ for \ }}

\def\io{\int_\Omega}
\def\iop{\int_{\Omega'}}

\def\dive{\mbox{div\,}}

\def\dd{\,\mbox{\rm d}}

\def\bchi{\chi}

\def\parchi{\frac{\partial}{\partial\bchi}}

\def\pardt{\frac{\partial}{\partial t}}

\def\vp{\varphi}

\def\matr#1#2#3#4{\left(\!\!\begin{array}{cc} #1 & #2\\#3 & #4\end{array}\!\!\right)}
\def\vect#1#2{\left(\!\!\begin{array}{c} #1\\#2\end{array}\!\!\right)}

\def\be{\begin{equation}\label}
\def\ee{\end{equation}}

\def\barr{\begin{array}}
\def\earr{\end{array}}

\def\ber{\begin{eqnarray}}
\def\eer{\end{eqnarray}}

\def\bers{\begin{eqnarray*}}
\def\eers{\end{eqnarray*}}

\numberwithin{equation}{section}

\begin{document}

\title{\bf A nonlocal \red{quasilinear} multi-phase system\\
 with nonconstant specific heat\\
and heat conductivity\footnote{\red{{\bf Acknowledgment.} The authors gratefully
acknowledge financial support by the MIUR-PRIN Grant
2008ZKHAHN ``Phase transitions, hysteresis and multiscaling'',
the GA\v{C}R Grant P201/10/2315, the program SMART-MATH of CNR/AV\v{C}R,
the MIUR-PRIN Grant 20089PWTPS ``Mathematical Analysis for inverse problems towards applications'',
the DFG Research Center {\sc Matheon} in Berlin, and the IMATI of CNR in Pavia.}}}

\author{\renewcommand{\thefootnote}{\! $\fnsymbol{footnote}$}
Pierluigi Colli\footnote{Dipartimento di Matematica, Universit\`a di Pavia,
Via Ferrata, 1, 27100 Pavia, Italy, E-mail {\tt pierluigi.colli@unipv.it}}, \
Pavel Krej\v{c}\'{\i}\footnote{Institute of Mathematics,
Czech Academy of Sciences, \v{Z}itn\'a 25, CZ-11567 Praha 1, Czech Republic,
E-mail {\tt  krejci@math.cas.cz}}, \
Elisabetta Rocca\footnote{Dipartimento di Matematica, Universit\`a di Milano,
Via Saldini 50, 20133 Milano, Italy, E-mail {\tt elisabetta.rocca@unimi.it}}, \ 
J\"urgen Sprekels\footnote{Weierstrass Institute for Applied
Analysis and Stochastics, Mohrenstr.~39, D-10117 Berlin,
Germany, E-mail {\tt  sprekels@wias-berlin.de}}}

\date{}

\maketitle

\renewcommand{\thefootnote}{}

\vspace{-.4cm}

\noindent {\bf Abstract.} In this paper, we prove the existence
and global boundedness from above for a solution to an
integrodifferential model for nonisothermal multi-phase
transitions under nonhomogeneous third type boundary conditions.
The \red{system} couples a \red{quasilinear} internal energy balance
ruling the evolution of the absolute temperature with a vectorial
integro-differential \red{inclusion} governing the (vectorial)
phase-parameter dynamics. The specific heat and the heat
conductivity $k$ are allowed to depend both on the order parameter
$\chi$ and on the absolute temperature $\teta$ of the system, and
the convex component of the free energy may or may not be
singular. \red{Uniqueness} and continuous data dependence are
also proved \red{under additional assumptions}.
 \vspace{.4cm}

\noindent {\bf Key words:}\thspace Phase transitions, nonlocal
models, \red{quasilinear} integrodifferential vectorial
equation\vspace{4mm}

\noindent {\bf AMS (MOS) subject clas\-si\-fi\-ca\-tion:}\thspace
\red{35K51, 35K59, 35K65, 45K05, 80A22}

\pagestyle{empty}



\pagestyle{myheadings}
\newcommand\testopari{\sc }
\newcommand\testodispari{\sc Colli-Krej\v{c}\'{\i}-Rocca-Sprekels}
\markboth{\testodispari}{\testopari}

\bigskip

\section{Introduction}
\label{intro}

In this paper we consider a nonisothermal multi-phase transition process occurring
in a bounded container $\Omega\subset \RR^N$, $N\in \NN$, with Lipschitz boundary $\partial\Omega$.
The state variables describing the evolution of the system are the absolute temperature
$\teta>0$ and the vectorial order parameter $\chi\in \RR^d$, $d\in \NN$.
Following the idea that was already described in the pioneering papers
\cite{vdW} and \cite{cahnhill}, but which has been only recently
analyzed in a more systematic way (cf., e.\,g.,
\cite{bc}--\cite{bates}, \cite{ChF}--\cite{ckrs}, \cite{fip}--\cite{ksPF},
\cite{sz}), we take into account long range interactions between particles.
Then \red{the model equations resulting from the energy and entropy balance relations} have the form
\begin{align}\label{eq1}
&(e(\teta,\chi))_t
+\left(\lambda(\chi)+\beta\vp(\chi)\right)_t+b[\chi]\chi_t
-\dive{\left(k(\teta,\chi)\nabla\teta\right)}=0\quad\hbox{in $Q_{\infty}:=\Omega\times(0,+\infty)$}, \\
\label{eq2}
&\mu(\teta)\chi_t+\lambda'(\chi)+b[\chi]+(\beta+\teta)\dvp(\chi)+\teta\sigma'(\chi)
+e_\chi(\teta,\chi)-\teta s_\chi(\teta,\chi)\ni 0\quad\hbox{in $Q_\infty$}, \\
\label{bou}
&k(\teta,\chi)\nabla\teta\cdot{\bf n}+\gamma(\teta-\teta_\Gamma)=0\quad\hbox{on $\Sigma_\infty:=\partial\Omega\times (0,\infty)$},\\
\label{ini} &\teta(\cdot, 0)=\teta_0,\quad \chi(\cdot,
0)=\chi_0\quad\hbox{in }\Omega,
\end{align}
where
${\bf n}$ denotes the outward normal vector to $\partial\Omega$, and \eqref{eq2} has to be understood as
an inclusion in $\RR^d$, \red{where}  \red{$\dvp$ is a possibly multivalued subdifferential of  a general proper, convex,
and lower semicontinuous  function  $\vp:\RR^d\to \RR\cup\{+\infty\}$}. \red{The physical meaning of the functions $e$, $s$, $\lambda$, $\sigma$
and of the positive constant $\beta$
is explained in (\ref{dense})--(\ref{entr}),} while $b[\chi]$ (whose explicit form will
be given) represents the nonlocal operator acting on $\chi$. With an abuse of notation we have
used the symbols $\lambda'$, $\sigma'$, $e_\chi$, $s_\chi$ for gradient vectors in $\RR^d$ and
omitted the scalar product symbol between $\RR^d$ vectors (like $b[\chi]$ and $\chi_t$ in \eqref{eq1})
in order not to overburden the presentation.  The function $\mu$
in \eqref{eq2} represents the (bounded away from 0) mobility of the system, while $\gamma$ denotes \red{the heat transfer
coefficient through the boundary} $\partial\Omega$.  \red{The external temperature} $\teta_\Gamma$ is a sufficiently regular
boundary datum on $\Sigma_\infty$, and $\teta_0$, $\chi_0$
are supposed to be two given initial configurations.

The main novelty here
is to consider a multi-phase nonlocal phase field system
in the case when the specific heat \red{$c_V(\teta,\chi)=\partial_\teta e(\teta,\chi)$
and the heat conductivity $k(\teta,\chi)$}
are not constant and depend on both the variables $\teta$ and $\chi$.
Suitable regularity and growth conditions will be specified in the following section.

Let us only note that many typical expressions for $c_V$ in a two-phase system (i.e.~in case $d=1$)
can be  included in our analysis.
In the solid-liquid system mentioned above, for example,
we may have different values $c_V^0(\teta)$ in the solid and $c_V^1(\teta)$ in the liquid phase,
hence,
we may define $c_V(\teta,\chi) = c_V^0(\teta) +\chi(c_V^1(\teta)-c_V^0(\teta))$
(cf.~\cite[Section IV.4]{v}). The value of $\chi$ can be kept between
$0$ and $1$ \red{by setting $\vp = I_{[0,1]}$}
(the indicator function of $[0,1]$). The physically meaningful case in which
the behaviour of $c_V^0$ and $c_V^1$ are powers of $\teta$ ($\sim \teta^\alpha$, $\alpha\geq 1$) near zero and
bounded functions for large $\teta$'s can be covered by our analysis. Regarding
the heat conductivity $k$,  typical expressions of the type $k(\teta,\chi)
=K_1(\teta)\chi+K_2(\teta)(1-\chi)$, in case of a two phase system \red{with $\chi\in[0,1]$},
for  quite general functions $K_1$ and $K_2$, are also allowed here.

The main goal of this paper is to study \red{the global existence} of solutions to system (\ref{eq1})--(\ref{ini}), coupling
a suitable variational formulation of the \red{semilinear}
parabolic partial differential equation \eqref{eq1} \red{for $\teta$} to the
 \red{integrodifferential inclusion \eqref{eq2}
for $\chi$. We also prove some uniform
in time} upper bound for the absolute temperature of the system (see Theorem~\ref{thexi} below).
Uniqueness of solutions
\red{is} obtained \red{under additional assumptions, in particular, in case that the heat
conductivity $k$ in (\ref{eq1}) depends only on $\teta$ and not
on $\chi$}.

\red{Before} entering into the mathematical discussion of the
problem, let us give a brief derivation of the system
(\ref{eq1})--(\ref{ini}), emphasizing, in particular, the differences
between local and nonlocal models.

We assume here that the multi-phase transition process can be completely described by the
evolution of the state variables $\teta(x,t)>0$, which represents the absolute temperature
of the system, and the order parameter $\chi(x,t)$, which here is a vector in $\RR^d$.
We fix some constant reference temperature $\teta_c$, which will be assumed to be equal to 1, for
simplicity.

Inspired by the nonlocal Cahn-Hilliard
model studied by Gajewski in \cite{g}, we consider the following
nonlocal specific free energy
\begin{equation}\no
F[\teta,\chi]=f_0(\teta,\chi)+B[\chi],
\end{equation}
where $B$ is a potential that accounts for long range interaction between particles. More specifically,
given a bounded, symmetric kernel $\kappa\,:\,\Omega\times\Omega\to \RR$
and an even smooth function $G\,:\,\RR^d\to\RR$, we choose
\begin{equation}\label{bigB}
B[\chi](x,t)\,:=\,\io \kappa(x,y)\,G(\chi(x,t)-\chi(y,t))\,\dd y.
\end{equation}
Note that the local potential
$(\nu/2)|\nabla\chi|^2$ used often in the literature, see \cite{v} and the references therein,
can be obtained as a formal
limit as $n\to \infty$ from the nonlocal one with the choice
$G(\eta) = |\eta|^2/2$, $\kappa(x,y)
= n^{N+2} \widetilde \kappa(|n(x-y)|^2)$, where $\widetilde\kappa$ is a nonnegative function with
support in $[0,1]$
and $\nu = 1/N\int_{\RR^N} \widetilde \kappa(|z|^2)|z|^2\,\dd z$. This follows from the formula
\begin{align}\no
\int_\Omega n^{N+2} \widetilde\kappa(|n(x-y)|^2)\,|\chi(x) - \chi(y)|^2\,\dd y
&= \int_{\Omega_n(x)} \widetilde\kappa(|z|^2) \left|\frac{\chi\left(x+ \frac{z}{n}\right)
- \chi(x)}{\frac{1}{n}}\right|^2\,\dd z\\ \no
\stackrel{n\to\infty}{\longrightarrow}&
\int_{\RR^N} \widetilde\kappa(|z|^2)\duavg{\nabla \chi(x), z}^2\,\dd z
= \nu |\nabla \chi(x)|^2
\end{align}
for a sufficiently regular $\chi$, where we denote
$\Omega_n(x) =n(\Omega - x)$. We have used the identity
$\int_{\RR^N} \widetilde\kappa(|z|^2)\duavg{v,z}^2\,\dd z = 1/N\,\int_{\RR^N}
\widetilde\kappa(|z|^2)|z|^2\,\dd z$ for every unit vector $v \in \real^N$
(cf.~the Introduction of \cite{krs2} for further details on this topic).

\red{Let $E$ and $S$ be the total energy and entropy densities, respectively.} The process is governed by the internal energy
and entropy balance relations
over an arbitrary control volume $\Omega'\subset\Omega$,
\begin{equation}\label{intbalg}
\frac{\dd}{\dd t}\iop E(\teta,\chi) \,\dd x+ \iopd\duavg{{\bf q}, {\bf n}} \,\dd s(x)
 = \Psi(\Omega')\,,
\end{equation}
\begin{equation}\label{claudu}
\frac{\dd}{\dd t}\iop S(\teta,\chi)\,\dd x +  \iopd\duavg{\frac{{\bf q}}{\teta}, {\bf n}}\, \dd s(x)\ge \ 0\,,
\end{equation}
where ${\bf q}$ is the heat flux vector, ${\bf n}$ is the unit outward
normal to $\partial\Omega'$, and $\Psi(\Omega')$ is the energy exchange
through the boundary of $\Omega'$ due to the nonlocal interactions.
Since $B[\chi]$ is a potential field, it does not contribute
to the entropy production in the Clausius-Duhem
inequality (\ref{claudu}).

The local form of the entropy balance reads
\begin{equation}\no
\teta S_t(\teta,\chi)+\div {\bf q}-\frac{\duav{{\bf q},\nabla\teta}}{\teta}\geq 0,
\end{equation}
and it has to be understood in the regularity context of Theorem~\ref{thexi} \red{below}.
This is certainly satisfied if
\begin{align}\no
&\duav{\mbox{\bf q}, \nabla\teta} \leq0\,,\\
\no
&\teta S_t(\teta,\chi) + \dive \mbox{\bf q}\ge 0\,.
\end{align}
Assuming $\teta>0$ and a suitable regularity with respect
to time (this will have to be justified in the
next sections), we obtain from (\ref{intbalg}) that
\begin{equation}\label{etst}
\iop (E_t - \teta S_t) \,\dd x \le \Psi(\Omega')\,.
\end{equation}
Differentiating the identities $F = E - \teta S = f_0 + B[\chi]$ with respect
to $t$, we obtain
\begin{equation}\label{freeendif}
F_t = E_t - \teta S_t - \theta_t S = \partial_\theta f_0 \theta_t
+\partial_\chi f_0 \chi_t + B[\chi]_t\,,
\end{equation}
where $\partial_\chi f_0$ stands for an element of Clarke's partial
subdifferential of $f_0$ with respect to $\chi \in \RR^d$, and
$\partial_\theta f_0$ is the partial derivative of $f_0$ with respect to
$\theta\in \RR$. Consequently,
\begin{equation}\label{thermo}
S = - \partial_\theta f_0 = s_0\,, \quad E = e_0 + B[\chi]\,, \quad
f_0 = e_0 - \theta s_0\,,
\end{equation}
and inequality (\ref{etst}) reads
\begin{equation}\label{etst2}
\iop (\partial_\chi f_0 \chi_t + B[\chi]_t) \,\dd x \le \Psi(\Omega')\,.
\end{equation}
The nonlocal interaction takes place only inside the domain $\Omega$,
hence $\Psi(\Omega) = 0$. A canonical way to satisfy these conditions
independently of the evolution of $\chi$ consists in choosing
the order parameter dynamics in the form
\begin{equation}\label{order}
\mu(\teta)\chi_t\in -\red{D_\chi}{\cal F}[\teta,\chi]
\end{equation}
with a factor $\mu(\theta) > 0$, where we denote
\begin{equation}\no
{\cal F}[\teta,\chi]\,=\,\io \, F[\teta,\chi]\,\dd x
\end{equation}
and  $\red{D_{\chi}}{\cal F}$ stands for the Clarke subdifferential of
${\cal F}$ with respect to the variable $\chi \in L^2(\Omega; \RR^d)$.
The inclusion sign in
(\ref{order}) accounts for the fact that $f_0(\teta,\chi)$ \red{includes}
terms that \red{are possibly not} Fr\'echet differentiable.
Condition (\ref{order}) is based on the assumption that the
system tends to move towards local minima of the free energy with a speed
proportional to $1/\mu(\teta)$.
Denoting
\begin{equation}\label{smallB}
b[\chi](x,t)\,:=\,2\io \kappa(x,y)\,G'\left(\chi(x,t)-\chi(y,t)\right)\,\dd y\,,
\end{equation}
where again, with an abuse of notation, $G'$ \red{stands for} the $d$-component vector $\nabla G$, we see that
the inequality (\ref{etst2}) holds
without prescribing any relationship between $\mu(\teta)$ and $B[\chi]$,
provided that we choose $\Psi(\Omega')$ in (\ref{intbalg}) as
\begin{equation}\label{psi}
\Psi(\Omega') = \iop \left(-b[\chi]\chi_t+B[\chi]_t\right) \dd x\,.
\end{equation}
The differential form of the energy balance (\ref{intbalg}) then reads
\begin{equation}\label{intenbal}
E_t+\dive{\bf q}\ =\ -b[\chi]\chi_t+B[\chi]_t\,.
\end{equation}

The specific heat $c_V(\theta,\chi)$ is the only thermodynamic state function, which can be identified from the measurements, \red{while the} local
internal energy and entropy densities are computed from the formulas
\begin{align}
\label{dense}
&e_0(\teta,\chi)=e_0(0,\chi)+e(\teta,\chi),\quad e(\teta,\chi)=\int_0^\teta c_V(\tau,\chi)\,\dd\tau\,,\\
\label{denss1}
&s_0(\teta,\chi)=s_0(0,\chi)+s(\teta,\chi), \quad s(\teta,\chi)=\int_0^\teta\frac{c_V(\tau,\chi)}{\tau}\,\dd\tau\,,
\end{align}
where $e_0(0,\chi), s_0(0,\chi)$ are in fact ``integration constants'',
which we choose as
\begin{align}\label{inten}
&e_0(0,\chi)=\lambda(\chi)+\beta\varphi(\chi),\\
\label{entr}
&s_0(0,\chi)=-\sigma(\chi)-\vp(\chi)\,,
\end{align}
where \red{$\vp:\RR^d\to\RR\cup\{+\infty\}$ is proper, convex,
and lower semicontinuous},
the functions $\lambda$ and $\sigma$ \red{are  sufficiently regular} on $\domvp$, and
the parameter $\beta$ is a positive constant.

Then, the free energy functional $F$ has the form
\begin{equation}\label{freeen}
F[\teta,\chi]\ = \ f(\teta,\chi)+
\lambda(\chi)+B[\chi]+(\beta+\teta)\varphi(\chi)+\teta\sigma(\chi)\,,
\end{equation}
where $f(\teta,\chi)=e(\teta,\chi)-\teta s(\teta,\chi)$.

Using \eqref{freeen}, \red{we rewrite} the phase dynamics
\eqref{order} as
\begin{equation}\label{phase}
\mu(\teta)\chi_t
+\lambda'(\chi)+b[\chi]+(\beta+\teta)\dvp(\chi)+\teta\sigma'(\chi)
+e_\chi(\teta,\chi)-\teta s_\chi(\teta,\chi)\ni 0\,,
\end{equation}
while the internal energy balance \eqref{intenbal} can be reformulated as
\begin{align}\label{teta}
&(e(\teta,\chi))_t
+\left(\lambda(\chi)+\beta\vp(\chi)\right)_t+b[\chi]\chi_t
-\dive{\left(k(\teta,\chi)\nabla\teta\right)}=0.
\end{align}

We now show that in the energy conserved case, that is, if we assume no-flux boundary conditions ($\gamma=0$ in \eqref{bou}),
the phase transition model with
a nonlocal interaction potential is compatible with the \"Ottinger-Grmela
GENERIC formalism (abbreviation for ``General Equation for the
Non-Equilibrium Reversible-Irreversible Coupling''), see \cite{og}.

Set
\ber\label{e19}
{\cal E}[\theta,\bchi](t) &=& \io E(\theta,\bchi)(x,t)\dd x\,,\\[2mm]\label{e20}
{\cal S}[\theta,\bchi](t) &=& \io S(\theta,\bchi)(x,t)\dd x\,,\\[2mm]\label{e20p}
{\cal B}[\bchi](t) &=& \io B[\bchi](x,t)\dd x\,.
\eer
We show that there exists a symmetric positive semidefinite matrix
${\bf M}[\theta,\bchi]$ such that
\be{e22}
{\bf M}[\theta,\bchi] \vect{D_\theta{\cal E}[\theta,\bchi]}{D_\bchi{\cal E}[\theta,\bchi]} =
\vect{0}{0}\,,
\ee
and such that the system (\ref{phase})--(\ref{teta}) has the form
\be{e21}
\frac{\partial}{\partial t} \vect{\theta}{\bchi}
= {\bf M}[\theta,\bchi] \vect{D_\theta{\cal S}[\theta,\bchi]}{D_\bchi{\cal S}[\theta,\bchi]}.
\ee
It suffices to choose (we omit the arguments of the state functions for simplicity)
\be{e23}
{\bf M}[\theta,\bchi] = \matr{M_0}{0}{0}{0} + \matr{m_{11}}{m_{12}}{m_{12}}{m_{22}},
\ee
where $M_0$ is the differential operator
\be{e24}
M_0 [y] = -\frac{1}{c_V}\dive\left(\theta^2 k(\theta,\bchi)\nabla
\frac{y}{c_V}\right)\,,
\ee
with homogeneous Neumann boundary condition, and $m_{ij}$ are scalars given by the formulas
\ber\label{e25}
m_{11} &=& \frac{\theta}{\mu(\theta)\,c_V^2} (D_\bchi{\cal E})^2\,,\\[2mm]\label{e26}
m_{12} &=& - \frac{\theta}{\mu(\theta)\,c_V} D_\bchi{\cal E}\,,\\[2mm]\label{e27}
m_{22} &=& \frac{\theta}{\mu(\theta)}\,.
\eer
Note that $m_{12}^2 = m_{11}m_{22}$ and $m_{11} \ge 0, m_{22}\ge 0$;
hence, ${\bf M}[\theta,\bchi]$ is positive semidefinite.
Furthermore, we have
\be{e28}
\vect{D_\theta{\cal E}}{D_\bchi{\cal E}} =
\vect{c_V}{D_\bchi{\cal E}} = \vect{c_V}{\parchi e_0+ D_\bchi{\cal B}}
, \quad \vect{D_\theta{\cal S}}{D_\bchi{\cal S}} =
\vect{c_V/\theta}{\frac{\partial S}{\partial\chi}}\,.
\ee
We easily check that (\ref{e22}) holds, and (\ref{e21}) has the form
\be{e29}
\vect{\theta_t}{\bchi_t}
= \vect{\dive (k(\theta,\bchi)\nabla\theta)/c_V}{0}
+ \matr{m_{11}}{m_{12}}{m_{12}}{m_{22}}\vect{c_V/\theta}{\frac{\partial S}{\partial\chi}}\,.
\ee
In component form, we have
\ber\label{e30}
\theta_t &=& \frac{1}{c_V}\left( \dive (k(\theta,\bchi)\nabla\theta)
+ \frac{1}{\mu(\theta)} D_\bchi{\cal E}\, D_\bchi{\cal F }\right),\\[2mm]\label{e31}
\bchi_t &=& -\frac{1}{\mu(\theta)} \, D_\bchi{\cal F}\,.
\eer
To see that (\ref{e30})--(\ref{e31}) coincides with (\ref{phase})--(\ref{teta}),
it suffices to take into account the formula
$$
\pardt e_0 = c_V\,\theta_t + \parchi e_0\,\chi_t = \dive (k (\theta,\bchi)\nabla\theta)
+ \frac{1}{\mu(\theta)} D_\bchi{\cal B}\, D_\bchi{\cal F}
= \dive (k (\theta,\bchi)\nabla\theta)
-  \chi_t\,D_\bchi{\cal B}\,.
$$
This shows that the model is compatible both with the standard principles of thermodynamics
and  the generalized thermodynamic formaliNotice thatsm introduced in \cite{og}. Note that Mielke
\cite{mielke} recently elaborated the GENERIC approach for the
phenomenology of thermoelastic dissipative materials.

\red{We} prove an existence result for a suitable variational formulation of
system (\ref{eq1})--(\ref{ini}). Using
a Moser technique, we also show that the temperature variable $\teta$
is globally bounded  from  above.
\red{The} uniqueness result holds true for particular classes of potentials $\vp$ provided that the  heat conductivity $k$ in
\eqref{eq1} does not depend on $\chi$, {\pier $\gamma=0$ in \eqref{bou},}  
and $c_V$ and $\mu$ satisfy a suitable growth condition around 0.

The paper is organized as follows. In Section~\ref{mainres}, we state our assumptions
on the data and our main results; in particular, global existence
for a suitable variational formulation of (\ref{eq1})--(\ref{ini}).
In Section~\ref{auxi}, we prove some auxiliary results related to the \blu{Lipschitz continuity
of solution operators to differential inclusions}. The proof will be developed as follows:
the problem is approximated by \red{partial time discretization}, regularization and a cut-off procedure
(cf. Subsection~\ref{approxi}).
Suitable a priori estimates (cf.~Subsection~\ref{apriori}) 
\red{allow us to pass to the limit with respect to the time step}
and regularization parameters,
while the cut-off is removed \red{by} proving an upper bound on the absolute temperature (which is independent of the
truncation parameter) by means of Moser techniques (cf.~Subsection~\ref{moser}). The uniqueness result is proved in Section~\ref{proofuni}.


\section{Main results}\label{mainres}

In this section, we state our main results on solvability conditions
for the system (\ref{eq1})--(\ref{ini}). We start by introducing a suitable
variational formulation; to this end, we
consider a bounded and Lipschitz domain $\Omega \subset \RR^N$, $N\ge 1$,
and for $t\in (0,\infty]$ we denote by $Q_t=\Omega\times (0,t)$ the open space-time cylinder
and by $\Sigma_t$ its lateral boundary $\partial\Omega\times (0,t)$.
We use, for the sake of simplicity, the same symbol $H$ for both
$L^2(\Omega)$ and $L^2(\Omega\,;\,\RR^N)$, while for an arbitrary integer $d$, $\mathbf{H}$ denotes
the space $L^2(\Omega;\RR^d)$. $H$ and $\mathbf{H}$  are both endowed with the standard scalar product
which we denote by $(\cdot,\cdot)$. The symbol  $V$ stands for the space $H^1(\Omega)$, and
$V'$ for its dual space, while the symbol $\mathbf{V}$ denotes the space $H^1(\Omega;\RR^d)$, $\duav{\cdot,\cdot}$ being
the duality $V'-V$ and \red{$\mathbf{V}'-\mathbf{V}$}.
Then, the following dense and continuous embeddings,
where we identify $H$ (and $\mathbf{H}$) with its dual space $H'$ (and $\mathbf{H}'$),
hold true: $V\hookrightarrow H\equiv H'\hookrightarrow V'$, and $\mathbf{V}\hookrightarrow
\mathbf{H}\equiv \mathbf{H}'\hookrightarrow \mathbf{V}'$.
Finally, we rewrite the system (\ref{eq1})--(\ref{ini}) in the following
variational formulation:
\begin{align}\label{p1}
&\duav{\partial_t(e(\teta,\chi)), z}+\io k(\teta,\chi)\nabla\teta\cdot\nabla z\dd x
+\int_{\partial\Omega}\gamma(\teta-\teta_\Gamma)\,z\dd A\\
\no
&\hspace{20mm}=-\io\left(\lambda'(\chi)\partial_t\chi
+\beta\left(\vp(\chi)\right)_t+
b[\chi]\chi_t\right)\, z\dd x\quad \forall z\in V,\quad\hbox{a.e.~in }(0,\infty)\,,\\
\label{p2}
&\mu(\teta)\chi_t+\lambda'(\chi)+\teta\sigma'(\chi)+
(\beta+\teta)\dvp(\chi)+b[\chi]\\
\no
&\hspace{24mm}+e_\chi(\teta,\chi)-\teta s_\chi(\teta,\chi) \ni 0\quad\hbox{a.e.~in }Q_\infty\,,
\end{align}
where \eqref{p2} has to be understood as an inclusion in $\RR^d$ with $b[\chi]$ defined by
\eqref{smallB}, and $e$ and $s$ are defined in (\ref{dense})--(\ref{denss1}).
Letting (cf.~\eqref{ini}) $u_0:=e(\teta_0,\chi_0)$,  we prescribe the initial conditions
\begin{align}
\label{p3}
&e(\teta,\chi)(0)=u_0,\quad\chi(0)=\chi_0\quad\hbox{a.\,e.~in }\Omega\,,
\end{align}
and suppose that the data fulfil the following assumptions.
\goodbreak

\behy{\rm (Existence)}
\label{hyp1}
Let us fix positive constants  $C_\sigma$, $C_\lambda$, $k_0$, $k_1$, $\underline{c}$, $\bar c$, $c_1$,
$\beta$,  $C_0$,
and assume that
\begin{itemize}
\item[{\rm (i)}] $\vp : \RR^d \to \RR\cup \{+\infty\}$ is a
proper, convex, and lower semicontinuous function, $\domvp$ is its domain;
\item[{\rm (ii)}] $\sigma,\,\lambda\in W^{2,\infty}(\domvp)$, $|\sigma'(r)|\leq C_\sigma$,
$|\lambda'(r)|\leq C_\lambda$ for all $r\in \domvp$;
\item[{\rm (iii)}] \red{$\kappa\in
W^{1,\infty}(\Omega\times\Omega)$,  $\kappa(x,y)=\kappa(y,x)$ a.\,e.~in $\Omega\times\Omega$, $G\in W^{2,\infty}(\domvp - \domvp)$,
$G(z)=G(-z)$ for all $z\in (\domvp - \domvp)$}\,;
\item[{\rm (iv)}] $k\,:\,\RR \times \domvp\to(0+\infty)$\hbox{ is a locally Lipschitz continuous function such that }
$0< k_0\leq k(v,w)\leq k_1$ for all $v\in \RR$ and $w\in \domvp$;
\item[{\rm (v)}] \blu{The function $\mu$ maps $[0,\infty)$ in $(0,\infty)$ and the function
$\theta \mapsto \frac{1+\theta}{\mu(\theta)}$} is \blu{bounded and} Lipschitz continuous
on $[0,\infty)$ \blu{with Lipschitz constant $L_\mu$};
\item[{\rm (vi)}] $c_V\,:\,[0,+\infty)\times\domvp\to[0,+\infty)$ is a continuous function satisfying
\begin{align}\label{c1}
&c_V(0,\chi)=0, \quad  0<c_V(\teta,\chi)\leq \bar c \quad\forall\teta\in (0,+\infty), \quad \forall \chi\in \domvp;\\
\label{c2}
&\underline{c}\leq c_V(\teta,\chi)\quad\forall(\teta,\chi)\in [1,+\infty)\times\domvp;\\
\label{c3}
&\hbox{the function } \teta\mapsto \frac{c_V(\teta,\chi)}{\teta}\quad\hbox{is integrable in $(0,1)$ for all } \chi\in\domvp.
\end{align}
Moreover, for all $(\teta,\chi)\in [0,+\infty)\times\domvp$ there exists the gradient  $(c_V)_\chi(\teta,\chi)$, and it holds
\begin{align}
\label{c4}
&|(c_V)_\chi(\teta,\chi)|\leq c_1 c_V(\teta, \chi),\quad|(c_V)_\chi(\teta,\chi_1)-(c_V)_\chi(\teta,\chi_2)|\leq c_1|\chi_1-\chi_2|\\
\no
&\mbox{ for all }\teta\in[0,+\infty), \,\, \chi, \,\chi_1,\,\chi_2\in \domvp.
\end{align}
\red{Let $e$ and $s$ be defined by formulas (\ref{dense})--(\ref{denss1})
and suppose that}
\begin{align}\label{s1}
&0<s(1,\chi)\leq c_1\quad\forall\chi\in\domvp;\\
\label{s2}
&|s_\chi(\teta_1,\chi_1)-s_\chi(\teta_2,\chi_2)|\leq c_1\left(|\teta_1-\teta_2|+|\chi_1-\chi_2|\right)\\
\no
&\mbox{ for all }\teta_1,\teta_2\in[0,+\infty), \,\, \chi_1,\,\chi_2\in \domvp.
\end{align}
\item[{\rm (vii)}] $\chi_0\in \mathbf{V}\cap L^\infty(\Omega)^d$. Moreover, for any $C>0$ set
$$\dcvp =
\{\chi\in \domvp\,:\, \exists\,\xi\in\dvp(\chi)\,:\,|\xi|\leq C\},$$ and assume that
$\chi_0(x)\in \dczvp$
a.\,e. in $\Omega$;
\item[{\rm (viii)}] $\teta_0,\,u_0\in L^\infty(\Omega)$ fulfil $u_0=e(\teta_0,\chi_0)$ and
$\teta_0(x)> 0$ a.\,e. in $\Omega$;
\item[{\rm (ix)}] $\gamma\in L^\infty(\partial\Omega)$ is a nonnegative
function;
\item[{\rm (x)}] $\teta_\Gamma\in L^\infty(\Sigma_\infty)$ is such that $\teta_\Gamma(x,t)>0$ a.e.
and $\log(\teta_\Gamma)\in L^1(\Sigma_\infty)$;
\item[{\rm (xi)}] Let \red{$\zeta\in (W^{1,1}_{loc}(0,\infty))^d$} be the solution to the differential
inclusion
\begin{equation}\label{diffinclu}
\alpha(t)\zeta_t+\dvp(\zeta)\ni g(t) \quad \hbox{\red{a.e.}},
\end{equation}
with the initial condition $\zeta(0)=\zeta_0$,
\red{$\zeta_0\in \RR^d$, and given data $g\in (L^\infty(0,\infty))^d$ and $\alpha\in L^\infty_{\rm loc}(0,\infty)$ such that
 $0<\alpha_0\leq \alpha(t)$  a.e. We assume that}
there exists a positive constant $D>0$ such that for all $C>0$ such that
$|g(t)|\leq C$,  and $\zeta_0\in \dcvp$, \red{we have}
\begin{equation}\label{diffesti}
|(g-\alpha\zeta_t)(t)|\leq DC \quad \hbox{\red{a.e.}}\,.
\end{equation}
\end{itemize}
\enhy

We are now in position to state the existence theorem.

\bete\label{thexi} \blu{\rm (Existence)}
Let Hypothesis~\ref{hyp1} hold.
Then there exists at least one pair $(\teta,\chi)$
that solves system (\ref{p1})--(\ref{p3}) and such that
\begin{align}\label{regteta}
&\teta\in L^\infty(Q_\infty)\cap L^2(0,\infty;V), \quad \red{(e(\teta,\chi))_t\in L^2_{loc}(0,\infty;V')}\,,\\
\label{positeta}
&\teta(x,t)> 0\quad\hbox{a.\,e.~in }Q_\infty\,,\\
\label{regchi}
& \chi\in L^\infty_{\rm loc}(Q_\infty)^d\grr{\cap L^\infty_{\rm loc}(0,\infty;\mathbf{V})},\quad
\chi_t\in L^\infty(Q_\infty)^d; \\
\no
&\exists\,C>0\,:
\,\chi(x,t)\in \dcvp\quad\hbox{a.\,e.~in }Q_\infty\,.
\end{align}
Moreover, there exists a positive constant $\overline\teta$
independent of $t$
such that the following uniform upper bound hold:
\begin{equation}\label{upperlower}
 \teta(x,t) < \overline\teta\quad\hbox{for a.\,e.~}
(x,t)\in Q_\infty.
\end{equation}
\ente

\red{
\behy{\rm (Uniqueness)}
\label{hypuni}
Assume that Hypothesis~\ref{hyp1} is satisfied and suppose moreover that
\item[{\rm (i)}] $k(\teta,\chi)=\bar k(\teta)$ for all $\teta\in\RR$ and $\chi\in \RR^d$;
\item[\grr{\rm (ii)}] Fix $T\in (0,\infty)$ and suppose that there exists a positive constant $R$, depending only on $C$, $\alpha_0$, and $T$
such that the solutions
$\zeta_1,\,\zeta_2\in W^{1,\infty}(0,T)$ to \eqref{diffinclu} associated with
data $\zeta_{01}, \zeta_{02}
\in \dcvp$, $\alpha_1, \alpha_2 \in L^\infty(0,T)$, and with
$g_1, g_2 \in L^\infty(0,T)$ complying with the constraint
\begin{equation}\label{hg}
|g_i(t)|\leq C,\quad i=1,\,2,\quad\hbox{a.\,e.~in }(0,T),
\end{equation}
satisfy for all $t\in(0,T)$ the inequality
\begin{equation}\label{cdg}
\itt|\dot\zeta_1-\dot\zeta_2|(\tau)\dd\tau + |\zeta_1-\zeta_2|(t)
\leq  R\Big(|\zeta_{01}-\zeta_{02}|+\itt\left( \left|\frac{1}{\alpha_1}-\frac{1}{\alpha_2}\right|(\tau)+|g_1-g_2|(\tau)\right)\dd\tau\Big)\,;
\end{equation}
\item[{\rm (iii)}] Define  $\tilde c(\teta):=\min\{c_V(v, \chi): \chi\in\domvp,\,\, v\geq \teta\}$ and
assume that $\displaystyle\int_0^1 \frac{\tilde c(v)\mu(v)}{v^2}\dd v=+\infty$;
\item[{\rm (iv)}] The function $v\mapsto v^2/\mu(v)$ is nondecreasing in $(0,+\infty)$;
\item[\grr{\rm (v)}] There exists $\teta_*>0$ such that $\teta_0(x)\geq \teta_* $ {\pier for a.e. $x\in\Omega$, and $\gamma\equiv 0$ on 
$\partial \Omega$};
\item[\grr{\rm (vi)}] \grr{Assume that}
\grr{$\tilde c(\teta) >0$ for every $\teta\in (0,\infty)$}.
\enhy
}

\beos\label{remcv}
\red{Hypothesis~\ref{hyp1} allows for physically meaningful choices of
$c_V$. We can choose, for example, a $(d+1)$-component model with $K=\{\chi_i\geq 0, \,\,\sum_{i=1}^d\chi_i\leq 1\}$, $\chi_0=1-\sum_{i=1}^d\chi_i$,
$\vp=I_K$ (indicator function of the set $K$), and $c_V(\teta,\chi)=\sum_{i=1}^d c_i(\teta)\chi_i$,
where $c_i(\teta)$ behave asymptotically at $0$ and $\infty$ like
$\displaystyle\frac{\teta^\alpha}{1+\teta^\alpha}$. Further examples of potentials $\vp$ complying with Hypothesis~\ref{hyp1}(xi) and Hypothesis~\ref{hypuni}(ii)
will be given in the following Section~\ref{auxi}.}
\eddos

We state then our last result regarding uniqueness and
continuous data dependence for (\ref{p1})--(\ref{p3}).

\bete\label{thuni} \blu{\rm (Uniqueness)}
Suppose that Hypothesis~\ref{hypuni} is satisfied.
Let $T\in (0,\infty)$ be fixed. Then, there exists a positive constant $\underline{\teta}(T)$ such that
\begin{equation}\label{lowbou}
\teta(x, t)\geq \underline{\teta}(T)\quad \hbox{for a.e. }(x, t)\in Q_T\,.
\end{equation}
Moreover, if $(\teta_1,\chi_1)$,
$(\teta_2,\chi_2)$ are two solutions to (\ref{p1})--(\ref{p3})
 in the sense of Theorem \ref{thexi}
associated with initial data $
\theta_{01},\chi_{01}$ and $
\theta_{02},\chi_{02}$, 
respectively, and  $\hat\theta = \theta_1- \theta_2$,
$\hat\chi = \chi_1-\chi_2$, $\hat\chi_0 = \chi_{01}-\chi_{02}$,
$\hat\theta_0 = \theta_{01}- \theta_{02}$, 
then, there exists
a constant $C_T>0$ such that
\begin{equation}\label{conti}
\int_0^T\io |\hat\theta(x,t)|^2\, \dd x\,\dd t +
\max_{t\in[0,T]} \io |\hat\chi(x,t)|^2\, \dd x
\le C_T \left(|\hat\theta_0|_H^2 + |\hat\chi_0|_{\bf H}^2
\right)\,.
\end{equation}
Finally, beside Hypothesis~\ref{hypuni}, assume that
\begin{equation}\label{regotetagamma}
{\pier \teta_0\in V }\,.
\end{equation}
Then, the $\teta$-component of the solution $(\teta,\chi)$ to (\ref{p1})--(\ref{p3}) has the further regularity
\begin{equation}\label{regotetat}
\teta\in L^\infty(0,T;V), \quad \teta_t\in L^2(0,T; H)\,.
\end{equation}
\ente


\section{A differential inclusion}
\label{auxi}

This section is devoted to the description of some properties of solutions to general differential
inclusions of the form \eqref{diffinclu}, which are used in the proof of Theorems~\ref{thexi}, \ref{thuni}.

\red{First we provide some examples
of  functions $\vp$ satisfying Hypothesis~\ref{hyp1}(xi),
and we prove some further properties for space and time dependent differential inclusions that follow
exactly from this assumption. Finally, we give examples of functions
$\vp$ satisfying Hypothesis~\ref{hypuni}(ii).}

\paragraph{Examples of functions complying with Hypothesis~\ref{hyp1}(xi).}

\bepr\label{exa} The function $\vp$ \red{introduced} in 
Hypothesis~\ref{hyp1}(i) also satisfies
Hypoth\-e\-sis~\ref{hyp1}(xi) \red{in each of} the following cases:
\begin{itemize}
\item[{\rm (a)}] if $d=1$;
\item[{\rm (b)}] if \red{$\vp$ is the indicator function $I_K$ associated with a closed, and convex set $K\subset\RR^d$.} In this case one has $D=1$;
\item[{\rm (c)}] if $\vp(x)=f(M_K(x))$, where \red{$f:[0,f_0)\to [0,+\infty)$} is an
increasing and convex $C^1$ function such that $f(0)=f'(0)=0$, \red{$f_0>0$},
and $M_K$ is the Minkowski functional of $K$, a closed, convex set in $\RR^d$ such that $B_r(0)\subset K\subset B_R(0)$, defined by the formula
$M_K(x)=\inf\left\{s>0\,;\,\frac1s x\in K\right\}$, for $x\in \RR^d$.
Then $D=R/r$.
\end{itemize}
\empr
The proof of point (a) follows directly from \cite[Prop.~3.4]{krs1}, (b) is obvious.
Let us prove the point (c). In order to do that, we first need to prove the following auxiliary result.
\bepr\label{Mink}
Let $\vp(x)$ be as in Proposition~\ref{exa}(c).
Then  for every $C>0$ there exists  $C_1>0$ such that for every $x\in \domvp$ we have the following implications
\begin{align}\label{implI}
&\vp(x)\leq C_1\,\Rightarrow\,\sup\{|\eta|\,:\,\eta\in\dvp(x)\}\leq (R/r)C,\\
\label{implII}
&\vp(x)\geq C_1\,\Rightarrow\,\inf\{|\eta|\,:\,\eta\in\dvp(x)\}\geq C.
\end{align}
\empr
\proof \red{We first prove the following equivalence
\begin{equation}\no
\eta\in\dvp(x)\quad\Leftrightarrow \quad\left(\eta=cw,\quad w\in\partial M_K(x),\quad c=f'(M_K(x))\right)\,.
\end{equation}
We clearly have $\dvp(0)=\{0\}$. For $x\neq 0$, take $\gamma\in (0,1)$. Then, for $\eta\in \dvp(x)$ and for all $y\in \domvp$, we have
$$\duav{\eta, x-(x-\gamma(x-y))}\geq \vp(x)-\vp(x-\gamma(x-y))\,,$$
hence
$$ \duav{\eta, x-y}\geq \frac{1}{\gamma}\left(f(M_K(y))-f(M_K(x))-\gamma (M_K(x)-M_K(y))\right)\,.$$
Letting $\gamma$ tend to $0$, we obtain that $\displaystyle\frac{\eta}{f'(M_K(x))}\in \partial M_K(x)$. Conversely,
for $w\in \partial M_K(x)$, we have
$$\duav{f'(M_K(x))w, x-y}\geq f'(M_K(x))(M_K(x)-M_K(y))\geq f(M_K(x))-f(M_K(y)),$$
which we wanted to prove.}

Let now $C$ be a given positive constant. For all $w\in \partial M_K(x)$ and  $x\neq 0$
we have $(1/R)\leq |w|\leq (1/r)$. From this we deduce
that if $f'(M_K(x))|w|<C$, then $f'(M_K(x))<CR$, and so $M_K(x)<(f')^{-1}(CR)$ and $\vp(x)<f((f')^{-1}(CR))$.
We can choose $C_1=f((f')^{-1}(CR))$, and \eqref{implII} is proved.
Suppose now that $f(M_K(x))\leq C_1$. Then we have $M_K(x)\leq (f')^{-1}(CR)$, hence $f'(M_K(x))\leq CR$ and
$f'(M_K(x))|w|\leq (R/r)C$, and \eqref{implI} is proved.\QED

We conclude the proof of Proposition~\ref{exa}(c) by proving the following Proposition~\ref{proMink}.

\bepr\label{proMink}
Let $\vp$ be as in Proposition~\ref{exa}(c).
Then Hypothesis~\ref{hyp1}(xi) is satisfied with $D=R/r$.
\empr
\proof
Consider some $\zeta$ satisfying inclusion \eqref{diffinclu}  with initial datum $\zeta_0$.
Then, the following equality holds true for all
$t\in \red{(0,\infty)}$:
\begin{equation}\no
|\alpha(t)\zeta_t|^2+|g(t)-\alpha(t)\zeta_t|^2+2\alpha(t)\vp(\zeta)_t=|g(t)|^2;
\end{equation}
hence, we immediately deduce that, for all
$t\in \red{(0,\infty)}$,
\begin{equation}\label{eqvpchi}
\vp(\zeta)_t=\frac{1}{2\alpha(t)}\left(|g(t)|^2-|g(t)-\alpha(t)\zeta_t|^2-|\alpha(t)\zeta_t|^2\right).
\end{equation}
\red{By assumption, we have $|g(t)|\leq C$ and $g(t)-\alpha(t)\zeta_t\in \dvp(\zeta)$. In view of \eqref{implII}, there is a constant $C_1$ such that
$$\vp(\zeta)_t\leq\frac{1}{2\alpha(t)}(C^2-C^2-|\alpha(t)\zeta_t|^2)\leq 0\quad \red{\forall t\in (0,\infty)}\quad\mbox{if}\quad\vp(\zeta)\geq C_1\,.$$}
Hence,
\begin{equation}\label{npvp}
\vp(\zeta)_t(\vp(\zeta)-C_1)^+\leq 0\quad\forall t\in (0,\infty).
\end{equation}
Integrating from $0$ to $t$ and using the assumption $\zeta_0\in\dcvp$, it
follows that  $\vp(\zeta)(t)\leq C_1$ for all $t\in (0,\infty)$. Then, using \eqref{implI}, together with the fact that
$|g(t)|\leq C$ for all $t\in (0,\infty)$, we finally obtain that
\begin{equation}\label{esgchit}
|g(t)-\alpha(t)\zeta_t|\leq \frac{R}{r}C
 \quad\hbox{and}\quad |\alpha(t)\zeta_t|\leq \left(\frac{R}{r}+1\right)C,
\end{equation}
which concludes the proof.\QED

\paragraph{Properties of the solution mapping under
Hypothesis~\ref{hyp1}(xi).}

\red{
\bepr\label{lip1}
Let us consider the solutions
$\zeta_1,\,\zeta_2\in \red{W^{1,\infty}(0,\infty)}$ to \eqref{diffinclu} associated with
data $\zeta_{01}, \zeta_{02}
\in \dcvp$, $\alpha_1, \alpha_2 \in \red{L^\infty_{\rm loc}(0,\infty)}$, and with
$g_1, g_2 \in \red{L^\infty(0,\infty)}$ complying with the constraint
\begin{equation}\no
|g_i(t)|\leq C,\quad i=1,\,2,\quad\hbox{a.\,e.},
\end{equation}
and let Hypothesis~\ref{hyp1}(xi) hold.
Then there exists a constant $L$ such that for every $t\in (0,\infty)$ we have
\begin{equation}\label{eslip1}
|\zeta_1-\zeta_2|(t)\leq |\zeta_{01}-\zeta_{02}|+L\itt \left(\left|\frac{1}{\alpha_1}-\frac{1}{\alpha_2}\right|
+|g_1-g_2|\right)\dd\tau\,.
\end{equation}
\empr
\proof Test the difference of the two inclusions \eqref{diffinclu} by $\zeta_1-\zeta_2$ and divide the resulting
inequality by $\alpha_1$. Then we obtain for a.e. $t\in (0,\infty)$:
$$\duav{\dot{\zeta}_1-\dot{\zeta}_2, \zeta_1-\zeta_2}\leq  \left|\frac{1}{\alpha_1}
-\frac{1}{\alpha_2}\right||\duav{\alpha_2\dot{\zeta}_2, \zeta_1-\zeta_2}|+|\duav{g_1-g_2,\zeta_1-\zeta_2}|\,. $$
Using the bound for $|\alpha_2\dot{\zeta}_2|$ (cf. Hypothesis~\ref{hyp1}(xi), \eqref{diffesti}), we get
$$ \frac{\dd}{\dd t} |\zeta_1-\zeta_2|\leq L \left(\left|\frac{1}{\alpha_1}-\frac{1}{\alpha_2}\right|
+|g_1-g_2|\right)$$
from which \eqref{eslip1} immediately follows by integrating over $(0,t)$.
\qed
\bepr\label{lip2}  Let Hypothesis~\ref{hyp1}(xi) hold, and let $\zeta_n$ and $\zeta$ be the solutions of \eqref{diffinclu}
corresponding to \red{the} data ($g_n$, $\alpha_n$, $\zeta_{0n}$) and ($g$, $\alpha$,  $\zeta_0$), respectively, with
$|g_n(t)|\leq C$, $\zeta_{0n}\in {\cal D}_C(\vp)$.
If $\{\zeta_{0n}\}$ converges to $\zeta_0$ in $\RR^d$,
$\{g_n\}$ converges to $g$ and $\{\alpha_n\}$ converges to $\alpha$
in $L^2(0,T)$ for some $T>0$, then $\{\dot{\zeta}_n\}$ converges strongly to $\dot{\zeta}$ in $L^2(0,T)$.
\empr
\proof Test \eqref{diffinclu}, written for $\zeta_n$, by $\dot{\zeta}_n$ \red{in order to obtain}
$$\left(g_n-\alpha_n\dot{\zeta}_n\right)\dot\zeta_n=\frac{\dd}{\dd t}\vp(\zeta_n)\,.$$
\red{Now set} $\displaystyle\eta_n=\frac{g_n}{\sqrt{\alpha_n}}-2\sqrt{\alpha_n}\dot\zeta_n$. Then, by straightforward computations, we obtain that
$$\left|\frac{g_n}{\sqrt{\alpha_n}}\right|^2-|\eta_n|^2=4\frac{\dd}{\dd t}\vp(\zeta_n)\,.$$
\red{We know, by Proposition~\ref{lip1}, that $\{\zeta_n\}$ converges uniformly to $\zeta$ and that $\vp$ is Lipschitz continuous on
$\dcvp$. Hence, integrating over $(0,t)$,} we obtain that $|\eta_n|_{L^2(0,T)}\to|\eta|_{L^2(0,T)}$.
Since we know that $\eta_n\to\eta$ weakly in $L^2(0,T)$
(since $\dot\zeta_n\to\dot\zeta$ weakly in $L^2(0,T)$),
we infer $\eta_n\to\eta$ strongly in $L^2(0,T)$, which is sufficient in order to conclude the desired convergence.
\QED
}

\paragraph{Examples of functions complying with Hypothesis~\ref{hypuni}(ii).}

\bepr\label{exb}
\red{The function $\vp$ \red{introduced} in 
Hypothesis~\ref{hyp1}(i) satisfies
Hypothesis~\ref{hypuni}(ii)
in each of the following cases:
\begin{itemize}
\item[{\rm (a)}] if $d=1$;
\item[{\rm (b)}] if, for any $C>0$, \red{${\vp}$ is a $C^1$-function with
Lipschitz continuous derivative on $\dcvp$};
\item[{\rm (c)}] if $\vp=I_K$, where $K$ is either a polyhedron  or \red{a smooth convex set with nonempty interior}.
\end{itemize}
}
\empr

\proof \red{The proofs of (a) and (c) follow respectively from  \cite[Prop.~3.4]{krs1} and \cite[Thm.~7.1, p.~88]{dkt}.
We briefly show \red{here} how \red{to} proceed to prove case (b). Let us consider solutions
$\zeta_1,\,\zeta_2\in W^{1,\infty}(0,T)$ to \eqref{diffinclu} associated with \red{the}
data $\zeta_{01}, \zeta_{02}
\in \dcvp$, $\alpha_1, \alpha_2 \in L^\infty(0,T)$, and with
$g_1, g_2 \in L^\infty(0,T)$ complying with the constraint
\begin{equation}\no
|g_i(t)|\leq C,\quad i=1,\,2,\quad\hbox{a.\,e.~in }(0,T).
\end{equation}
By Hypothesis~\ref{hyp1}(xi), $\zeta_1, \zeta_2$ remain in ${\cal D}_{DC}(\vp)$.
Using the Lipschitz continuity of $\vp'$ on ${\cal D}_{DC}(\vp)$, we obtain that there exists a positive constant
$Q$ such that the following inequality holds true a.e.:
$$\alpha_1|\dot{\zeta}_1-\dot{\zeta}_2|\leq |\alpha_1-\alpha_2||\dot{\zeta}_2|+Q|\zeta_1-\zeta_2|+|g_1-g_2|.$$
Dividing  by $\alpha_1$, and using the bound for $|\alpha_2\dot{\zeta}_2|$, from Hypothesis~\ref{hyp1}(xi) (cf. \eqref{diffesti}),
we get
\begin{align}\label{b1}
|\dot{\zeta}_1-\dot{\zeta}_2|&\leq(C+DC)\left|\frac{1}{\alpha_1}-\frac{1}{\alpha_2}\right| +Q|\zeta_1-\zeta_2|+|g_1-g_2|\\
\no
&\leq M\left(\left|\frac{1}{\alpha_1}-\frac{1}{\alpha_2}\right| +|\zeta_1-\zeta_2|+|g_1-g_2|\right)
\end{align}
for some positive constant $M$ (depending on $C, D, Q$). We can rewrite this inequality in the following convenient form, for $t\in (0,T)$,
$$\frac{\dd}{\dd t}\left({\rm e}^{-Mt}|\zeta_1-\zeta_2|\right)
\leq {\rm e}^{-Mt}\left(\left|\frac{1}{\alpha_1}-\frac{1}{\alpha_2}\right|
+|g_1-g_2|\right).$$
Integrating over $(0,t)$, and using the previous inequality \eqref{b1}, we deduce
$$ |\dot{\zeta}_1-\dot{\zeta}_2|(t)\leq L\left(|\zeta_{01}-\zeta_{02}|+\left|\frac{1}{\alpha_1}-\frac{1}{\alpha_2}\right|
+|g_1-g_2|+\itt \left(\left|\frac{1}{\alpha_1}-\frac{1}{\alpha_2}\right|
+|g_1-g_2|\right)\dd\tau\right)$$
for some positive constant $L$ (depending on $C, D, Q$). Integrating once more in time \red{we arrive at} the desired inequality \eqref{cdg}.}
\QED

\bepr\label{mink}
\red{Let $f:[0,f_0)\to [0,+\infty)$ be an
increasing, convex function with locally Lipschitz continuous derivative, 
$f(0)=f'(0)=0$,
and let $K$ be a closed, convex set of class $C^{1,1}$ such that 
$B_r(0)\subset K\subset B_R(0)$.
Then, $\vp(x)=f(M_K(x))$ has a Lipschitz continuous derivative on $\dcvp$ 
for any $C>0$, i.e., property~(b) in Proposition~\ref{exb} is satisfied.}
\empr

\proof
\red{Let $C>0$ be given. We denote ${\cal D}_C f=\{s \blu{\in (0,f_0)}: f'(s)\leq RC\}$, and let $L_C$ be the Lipschitz constant of $f'$ on ${\cal D}_Cf$.
For $x\in \dcvp$ we have $|\vp'(x)|\leq C$, hence $f'(M_K(x))\leq RC$, that is, $M_K(x)\in {\cal D}_Cf$. We now estimate the difference
$|\vp'(x)-\vp'(y)|$ on $\dcvp$. Assume first that $x\neq 0$, $y=0$. Then
\begin{align}
\no
&|\vp'(x)-\vp'(y)|=|\vp'(x)|=f'(M_K(x))|M_K'(x)|\leq \frac{1}{r} f'(M_K(x))\\
\no
&\leq \frac1rL_CM_K(x)\leq\frac{1}{r^2}L_C|x|=\frac{L_C}{r^2}|x-y|\,.
\end{align}
Consider now the case $x\neq 0$, $y\neq 0$ and set $J_K(x)=M_K(x)M_K'(x)$,  $J_K(y)=M_K(y)M_K'(y)$. The mapping $J_K$ is Lipschitz continuous on
$\RR^d$ (with Lipschitz constant $L_J$) (see \cite[Section 5.2]{dkt}), and we have
\begin{align}\no
&|\vp'(x)-\vp'(y)|=|f'(M_K(x))M_K'(x)-f'(M_K(y))M_K'(y)|\leq \frac{f'(M_K(x))}{M_K(x)}|J_K(x)-J_K(y)|\\
\no
&\quad+|M_K'(y)|\frac{f'(M_K(x))}{M_K(x)}|M_K(x)-M_K(y)|+|M_K'(y)||f'(M_K(x))-f'(M_K(y))|\\
\no
&\leq L_C\left(L_J+\frac{2}{r^2}\right)|x-y|\,,
\end{align}
from which  the assertion follows.}
\QED

\red{A relevant case for applications is, for example, $\vp(x)=-\log(1-M_K^2(x))$, see \cite{gz}.}


\section{Existence of solutions}
\label{proofexi}

This section is devoted to the proof of the existence
result stated in Section~\ref{mainres}. We use a
technique based on approximations, a priori estimates,
and passage to the limit.

Let us first write down our equations (\ref{p1}) and (\ref{p2}) as
\begin{align}\label{p1new}
&\duav{ (e(\teta, \chi))_t, z}+\io k(\teta,\chi)\nabla\teta\cdot\nabla z\dd x
+\int_{\partial\Omega}\gamma(\teta-\teta_\Gamma)\,z\dd A\\
\no
&=-\io\left(\lambda'(\chi)\chi_t
+\beta\left(\vp(\chi)\right)_t+b[\chi]\chi_t\right)\,z\dd x\quad\forall z\in V,\quad\hbox{a.e.~in }(0,\infty),\\
\label{p2new}
&\mu(\teta)\chi_t
+(\beta+\teta)\dvp(\chi)\ni-\lambda'(\chi)-\teta\sigma'(\chi)-b[\chi]
-e_\chi(\teta,\chi)+\teta s_\chi(\teta,\chi)\\
\no
&\hspace{11.80cm}\hbox{a.e.~in }Q_\infty.
\end{align}


\subsection{\blu{A}pproximation}\label{approxi}

Assuming Hypothesis~\ref{hyp1} to hold, we proceed as follows: first
we extend \blu{the domain of definition of} $c_V(\teta,\chi)$ \blu{by putting}
$\widetilde{c}_V(\teta,\chi)=c_V(|\teta|,\chi)$ \blu{for $(\theta,\chi) \in \RR \times \domvp$, and set}
$$\widetilde{e}(\teta,\chi)=\int_{0}^\teta \widetilde{c}_V(\xi,\chi)\, \dd\xi\quad \blu{\hbox{ for }
(\theta,\chi) \in \RR \times \domvp\,.}
$$
\blu{We now fix a truncation parameter $\vr \ge 1$, which will be determined below, and define}\red{}
$$
\red{\widetilde{\mu}_\vr(\teta)=
\begin{cases}
\displaystyle \mu(|\teta|) &\quad\hbox{if }|\teta|\leq\vr\,\\
\displaystyle   \mu(\vr)(|\teta|-\vr)           &\quad\hbox{if }|\teta|\geq\vr
\end{cases}\,,}
$$
$$s_\chi^\vr(\teta,\chi)=
\begin{cases}
\displaystyle\int_{0}^{\teta}\frac{(\widetilde{c}_V)_\chi(\xi,\chi)}{\xi}\, \dd\xi&\quad\hbox{if }|\teta|\leq\vr\,\\
\displaystyle\int_{0}^{\vr}\frac{(\widetilde{c}_V)_\chi(\xi,\chi)}{\xi}\, \dd\xi&\quad\hbox{if }|\teta|\geq\vr
\end{cases}\,.
$$
\blu{We fix an arbitrary $T>0$, and split the interval $[0,T]$ into an equidistant partition
$0 = t_0, t_1, \dots, t_n$, $t_j = jT/n$ for $j=0,1,\dots, n$, $n\in \NN$, with the intention to let $n$ tend to $\infty$.
We choose sequences $\{\teta_{0,n}\}_n\in V$ and $\{\teta_{\Gamma, n}\}_n \in W^{1,2}(0,T; L^2(\partial\Omega))$ of approximate data
such that $\teta_{0,n}(x)\geq 1/n$ a.e.,  $\teta_{\Gamma,n}(x,t)\geq 1/n$ a.e., $\teta_{0,n} \to \teta_0$ strongly in $H$,
and $\teta_{\Gamma, n}\to \teta_{\Gamma}$ strongly
in $L^{2}(0,T; L^2(\partial\Omega))$.}

\blu{An approximate solution $(\theta_n, \chi_n)$ will be constructed successively in intervals $[t_{j-1}, t_j]$
for $j=1, \dots, n$. Assuming that it is already known on $[0, t_{j-1}]$, we define
\begin{equation}\label{bartetachi}
\bar\chi_n(x,t) = \chi_n(x, t_{j-1}) \quad x\in \Omega, \quad t\in (t_{j-1}, t_j),\quad  j=1,\dots, n \,,
\end{equation}
$$
\bar\teta_n(x,t)=
\begin{cases}
\teta_{0,n}&\quad\hbox{for } t\in [0,t_1)\\
\displaystyle\frac{n}{T}\int_{t_{j-2}}^{t_{j-1}} \teta_n(x, \tau)\, \dd \tau &\quad\hbox{for }
 t\in [t_{j-1}, t_j), \quad j\geq 2\,.
\end{cases}
$$
With this notation, we then state} the following approximating problem.
We use only the index $n$ for the variables here (\blu{omitting} the $\vr$ dependence),
for simplicity.
\smallskip

\noindent
{\sc Problem {(P)$_{(n, \vr)}$.}} \blu{Find} two functions $\teta_n\in H^1(0,T;H)\cap L^\infty(0,T;V)$ and
$\chi_n\in L^{\infty}(\Omega\times(0,T))^d$, $\dt\chi_n\in L^{\infty}(\Omega\times(0,T))^d$,
 such that $\teta_n\geq \e_n$ a.e. in $Q_T$ for some $\e_n>0$, $\chi_n\in {\cal D}_C(\vp)$,  and
for all $t\in (0,T)$  and $z\in V$, we have
\begin{align}\label{p1tau}
&\io\dt\left(\red{\frac{1}{n}\teta_n(t)}+\widetilde{e}(\teta_n(t), \chi_n(t))\right)\, z\, \dd x
+\io \red{k(\bar\teta_n(t),\bar \chi_n(t))}\nabla\teta_n(t)\cdot\nabla z\dd x
\\
\no
& +\int_{\partial\Omega}\gamma(\teta_n(t)- \blu{\teta_{\Gamma, n}(t)})z\dd A=
-\io\left((\lambda'(\chi_n)(t)+b[\chi_n](t))\dt\chi_n(t)
+\beta\dt\left(\vp(\chi_n(t))\right)\right)z\dd x\,,\\
\label{p2tau}
& \red{\widetilde{\mu}_\vr}(\teta_n(t))\dt\chi_n(t)
+(\beta+|\teta_n(t)|)\dvp(\chi_n)(t)\ni
-\lambda'(\chi_n)(t)-|\teta_n(t)|\sigma'(\chi_n)(t)-b[\chi_n](t)\\
\no
& -\widetilde{e}_\chi(\teta_n(t),\chi_n(t))
+|\teta_n(t)| s_\chi^\vr(\teta_n(t),\chi_n(t))\quad\hbox{a.e.~in }\Omega\,,
\end{align}
\blu{with initial} conditions
\begin{equation}
\label{initau}
\displaystyle \teta_n(0)=\teta_{0,n},\quad \chi_n(0)=\chi_0\,.
\end{equation}

\bele\label{exiptau}
Under Hypothesis~\ref{hyp1}, \red{for each $\vr>0$ and $n\in \NN$ the}
{\sc Problem {(P)$_{(n,\vr)}$}} has a
unique solution $(\teta_n, \chi_n)$ with the \grr{required properties}.
\enle

\proof \red{On each interval $(t_{j-1}, t_j)$, we can proceed as in the proof of \cite[Thm.~2.2, p.~290]{krs2}.
We test a Galerkin approximation of \eqref{p1tau} by the approximation of $\dt \teta_n$.
The estimates are sufficient to pass to the limit in the Galerkin scheme and to obtain a solution
on each interval $(t_{j-1}, t_j)$. We only have to check that the initial conditions at
$t_1, t_2, \dots$ are well defined. Indeed, since on each interval $(t_{j-1}, t_j)$
we have $\teta_n\in H^1(t_{j-1}, t_j; H)\cap L^\infty (t_{j-1}, t_j; V)$, we also obtain that $t\mapsto \teta_n(t, \cdot)$ is weakly
continuous in $(t_{j-1}, t_j)$  for every $j$ with values in $V$. Moreover, $\chi_n$ is strongly continuous with
values in $L^\infty(\Omega)^d$, and there exists a positive constant $C$ (independent of $n$) such that
$\chi_n(t, \cdot)\in \dcvp$ on  $(t_{j-1}, t_j)$  for every $j=1,\dots, n$. Hence, we can define the initial conditions at $t=t_j$ by
$\teta_n(t_j)=\teta_n(t_j-)$.} \QED


\subsection{A priori estimates}\label{apriori}

In this subsection, we perform suitable a priori estimates
(independent of $n$) \red{for} the solution.
In the following, we will denote by $C$ any positive constant that depends only on the
data of the problem but may vary from line to line. \blu{In particular, it} will not
depend on the \blu{truncation} parameter $\vr$ and \blu{discretization parameter $n$. If
such a dependence takes place}, we use the symbol $C_\vr$ \blu{for a constant that depends on $\vr$,
but not on $n$}. Again,
the same symbols will denote constants \red{that may differ} from line to line.

Let us, for simplicity, \blu{in this subsection occasionally omit the indices $n$ and write simply
$\theta, \chi$ instead of $\theta_n, \chi_n$ if no confusion arises. We denote (note that $\teta_n > 0$} and so $\tilde e=e$)
\begin{equation}\label{defutau}
u_n(t):=\frac{1}{n}\teta_n(t)+{e}(\teta_n(t), \chi_n(t) )\quad\hbox{for } \blu{t\in (0,T)}\,.
\end{equation}

\paragraph{\blu{Estimate for $\chi_t$.}}
Equation \eqref{p2tau} is of the form \eqref{diffinclu} with
\begin{equation}\label{applaux1}
\alpha(t)=\tilde\alpha(\teta) \ := \
\frac{\widetilde{\mu}_\vr(\theta)}{\beta+\theta} \ \ge \
\frac{\mu_0(1+\theta)}{\beta+\theta} \ \ge \
\mu_0\,\min\left\{1,\frac{1}{\beta}\right\}\,,
\end{equation}
\begin{equation}\label{applaux2}
g(t)=\ell[\theta,\chi] := - \frac{1}{\beta+\theta}\left(
\theta\sigma'(\chi)+\lambda'(\chi)+b[\chi]+{e}_\chi(\teta,\chi)
-\teta s_\chi^\vr(\teta,\chi)\right).
\end{equation}
First, let us note that, using \eqref{c1}, \eqref{c4},  and \eqref{s1}, we get
\begin{align}\label{estisro}
s_\chi^\vr(\teta,\chi)&\leq\int_0^{\min\{\teta,\vr\}}\frac{|({c}_V)_\chi(\xi,\chi)|}{\xi}\,\dd\xi
\ \leq\
c_1\int_0^1\frac{c_V(\xi,\chi)}{\xi}\,\dd\xi+c_1\int_1^\vr\frac{\bar c}{\xi}\, \dd\xi\\[2mm]
\no
&\leq c_1 s(1,\chi)+c_1\bar c\log\vr
\ \leq\
c_1^2+c_1\bar c\log\vr\,.
\end{align}
Hence, owing to Hypothesis~\ref{hyp1}, \blu{we have}
\begin{align}\no
|\ell[\teta,\chi]|&\leq C_\sigma+\frac{1}{\beta}\left(C_\lambda+C_b\right)+\frac{\teta}{\beta+\teta}
\sup_{0\leq\xi\leq\teta}|({c}_V)_\chi(\xi,\chi)|+\frac{\teta}{\beta+\teta}|s_\chi^\vr(\teta,\chi)|\\
\no
&\leq C_\sigma+\frac{1}{\beta}\left(C_\lambda+C_b\right)+c_1\sup_{0\leq\xi\leq\teta}|c_V(\xi,\chi)|
+c_1^2+c_1\bar c\log\vr\\
\no
&\leq C_\sigma+\frac{1}{\beta}\left(C_\lambda+C_b\right)+c_1\bar c+c_1^2+c_1\bar c\log\vr\,,
\end{align}
where $C_b$ denotes here the upper bound for the operator $b$ defined in \eqref{smallB}.
Let us set
\begin{equation}\label{Cell}
C_{\ell,\vr}:=C_\sigma+\frac{1}{\beta}\left(C_\lambda+C_b\right)+c_1\bar c+c_1^2+c_1\bar c\log\vr\,.
\end{equation}
Then the conditions of Hypothesis~\ref{hyp1}(xi) are satisfied with the choice
$C=\max\{C_{\ell,\vr},\,C_0\}$, where $C_0$ is defined in 
Hypothesis~\ref{hyp1}(vii) and $C_{\ell,\vr}$ is defined in \eqref{Cell}.
Hence,  we obtain the following estimates on $\chi$:
\begin{equation}\label{estichi}
|\chi_n|_{L^\infty(Q_T)}+|\dt\chi_n|_{L^\infty(Q_T)}+
|\dt(\varphi(\chi_n))|_{L^\infty(Q_T)} \ \le \ C(1+\log\vr)^2\,,
\end{equation}
where now $C$ is a constant independent of $\vr$.

\paragraph{\blu{Estimate for $\theta$.}}
Taking $z=\teta_n$ in \eqref{p1tau}, we get
\begin{align}\label{e1}
&\io\partial_t\left(\frac{1}{n}\teta(t)+{e}(\teta(t), \chi(t))\right)\,\teta(t)\dd x +\io
k(\bar\teta_n(t),\bar\chi_n(t))\nabla\teta(t)\cdot\nabla\teta(t)\dd x
\\
\no
&+\int_{\partial\Omega}\gamma(\teta(t)-\teta_{\Gamma, n}(t))\teta(t)\dd A\\
\no
&=
-\io\left(\lambda'(\chi)(t)\dt\chi(t)
+\beta\dt\left(\vp(\chi(t))\right)+b[\chi](t)\dt\chi(t)\right)\teta(t)\dd x\,.
\end{align}
Define
$$
U(\teta,\chi)=\int_0^\teta c_V(v, \chi)v \, \dd v \quad \blu{ \hbox{for } (\theta, \chi) \in (0,\infty)\times \domvp}\,.
$$
We have $\partial_t U(\theta, \chi) = \theta \partial_t e(\theta, \chi) + (\partial_\chi U - \theta \partial_\chi e)\chi_t$.
We integrate \eqref{e1} from $0$ to $t$ and rewrite the first term as follows:
\begin{align}\no
&\itt \io\partial_t\left(\frac{1}{n}\teta(\tau)+{e}(\teta(\tau), \chi(\tau))\right)\,\teta(\tau)\dd \tau\dd x \\
\no
&= \blu{\io\left(\frac{1}{2n}\teta^2(t) + U(\theta(t), \chi(t))\right)\dd x -
\io\left(\frac{1}{2n}\teta^2(0) + U(\theta(0), \chi(0))\right)\dd x}\\
\no
& \ \ - \itt\io (\partial_\chi U(\teta(\tau), \chi(\tau)) - \theta(\tau)
\partial_\chi e(\teta(\tau), \chi(\tau)))\chi_t(\tau)\dd \tau \dd x\,.
\end{align}
There exist two constants $C_1, C_2$ such that $U(\teta,\chi)\geq C_1\teta^2-C_2$.
\blu{Hence, by \eqref{estichi} and the Gronwall's lemma, we obtain}
\begin{equation}\label{e1fin1}
\|\teta_n\|_{L^2(0,T;V)\cap L^\infty(0,T;H)}\leq C(1+\log\vr)^2\,.
\end{equation}
By comparison, we also deduce that
\begin{align}
\label{e1fin2}
&\|\dt u_n\|_{L^2(0,T; V')}\leq C(1+\log\vr)^2.
\end{align}

\paragraph{\blu{Estimate for $\nabla\chi$.}}
\blu{The function}
$$\displaystyle\teta\mapsto\frac{\beta+\teta}{\widetilde{\mu}_\vr(\teta)}$$
 is Lipschitz continuous in $\RR$
due to Hypothesis~\ref{hyp1}(v) and, with the help of the mean value theorem, 
it is straightforward to deduce that
\begin{align}\label{numeronuovo}
&|\ell[\teta_1,\chi_1]-\ell[\teta_2,\chi_2]|\leq |\sigma'(\chi_1)-\sigma'(\chi_2)|+\frac{C_\sigma}{\beta}|\teta_1-\teta_2|
+\frac{1}{\beta}|\lambda'(\chi_1)-\lambda'(\chi_2)|\\
\no
&\quad+\frac{C_\lambda}{\beta^2}|\teta_1-\teta_2|
+\frac{1}{\beta}|b[\chi_1]-b[\chi_2]|+\frac{C_b}{\beta^2}|\teta_1-\teta_2|+\frac{1}{\beta+\theta_1}\int_0^{\theta_1}c_1|\chi_1-\chi_2|\, \dd\xi\\
\no
&\quad
+\left|\frac{1}{\beta+\theta_1}\int_0^{\theta_1}(c_V)_\chi(\xi,\chi_2)\, \dd\xi
-\frac{1}{\beta+\theta_2}\int_0^{\theta_2}(c_V)_\chi(\xi,\chi_2)\, \dd\xi\right|\\
\no
&
\quad+\frac{1}{\beta}\left|s_\chi^\vr(\teta_1,\chi_1)-s_\chi^\vr(\teta_2,\chi_2)\right|
+\frac{|s_\chi^\vr(\teta_2,\chi_2)|}{\beta^2}|\teta_1-\teta_2|\\
\no
&\leq|\sigma'(\chi_1)-\sigma'(\chi_2)|+\frac{1}{\beta}\left(|\lambda'(\chi_1)-\lambda'(\chi_2)|+|b[\chi_1]-b[\chi_2]|\right)\\
\no
&\quad
+\frac{1}{\beta^2}\left(C_\sigma\beta+C_\lambda+C_b\right)|\teta_1-\teta_2|+\frac{c_1}{\beta}|\chi_1-\chi_2|
+\frac{1}{\beta}\left(c_1\bar c+c_1\bar c\right)|\teta_1-\teta_2|\\
\no
&
\quad+\frac{c_1}{\beta}\left(|\teta_1-\teta_2|+|\chi_1-\chi_2|\right)+\frac{1}{\beta^2}\left(c_1^2+c_1\log\vr\right)|\teta_1-\teta_2|\,.
\end{align}
Hence, we can apply estimate \eqref{eslip1} in Proposition~\ref{lip1}
to \eqref{p2tau} with (for $x,y\in \Omega$)
\begin{align}\no
&\alpha_1(t)=
\left(\frac{\widetilde{\mu}_\vr(\theta)}{\beta+\theta} \right)(x,t),\\
\no
&\alpha_2(t)=
\left(\frac{\widetilde{\mu}_\vr(\theta)}{\beta+\theta} \right)(y,t),\\
\no
&g_1(t)=\ell[\theta(x,t),\chi(x,t)] = -\frac{1}{\beta+\theta}\Big(
\theta\sigma'(\chi)+\lambda'(\chi)+b[\chi]+{e}_\chi(\teta,\chi)\\
\no
&\qquad\qquad\qquad\qquad\qquad\qquad-\theta s_\chi^\vr(\teta,\chi)\Big)(x,t),\\
\no
&g_2(t)=\ell[\theta(y,t),\chi(y,t)] = -\frac{1}{\beta+\theta}\Big(
\theta\sigma'(\chi)+\lambda'(\chi)+b[\chi]+{e}_\chi(\teta,\chi)
\\
\no
&\qquad\qquad\qquad\qquad\qquad\qquad-\theta s_\chi^\vr(\teta,\chi)\Big)(y,t),
\end{align}
and
\begin{equation}\no
\zeta_1=\chi(x, t),\quad \zeta_2=\chi(y,t),\quad x,y\in \Omega.
\end{equation}
Hence,  we obtain
\begin{align}\label{grad0}
|\chi(x,t)-\chi(y,t)|&\leq |\chi_0(x)-\chi_0(y)|+ \grr{\hat L_\vr\Big(}
\itt|\teta(x,s)-\teta(y,s)|\, \dd s\\
\no
&+ \itt\left(\left|b[\chi](x,s)-b[\chi](y,s)\right|+|\chi(x,s)-\chi(y,s)|\right)\, \dd s\Big)\,,
\end{align}
where \grr{$\hat L_\vr$ depends on $L$, $L_\mu$, and the constants on the right-hand side of
\eqref{numeronuovo}.}
Now, recalling \eqref{smallB}, we have, \blu{by \eqref{numeronuovo}, that}
\begin{align}\label{grad0bis}
|\chi(x,t)&-\chi(y,t)|\leq |\chi_0(x)-\chi_0(y)|+\grr{\hat L_\vr}\itt|\teta(x,s)-\teta(y,s)|\, \dd s\\
\no
&\quad+2\grr{\hat L}_\vr\itt\io\left|\kappa(x,z)\left(G'(\chi(x,s)-\chi(z,s)) -G'(\chi(y,s)-\chi(z,s))\right)\right|\, \dd z\, \dd s\\
\no
&\quad
+2 \grr{\hat L_\vr} \itt\io\left|G'(\chi(y,s)-\chi(z,s))\left(\kappa(x,z)-\kappa(y,z)\right)\right|\,\dd z\, \dd s\\
\no
&\quad +\grr{\hat L_\vr} \itt|\chi(x,s)-\chi(y,s)|\, \dd s\,.
\end{align}
Thus, in view of  Hypothesis~\ref{hyp1}(iii), we obtain that
\begin{align}\label{grad1}
|\chi(x,t)&-\chi(y,t)|\leq |\chi_0(x)-\chi_0(y)|+\grr{\hat L_\vr} \itt|\teta(x,s)-\teta(y,s)|\, \dd s\\
\no
&\quad
+ \grr{\hat L_\vr}(2L_b+1)\itt|\chi(x,s)-\chi(y,s)|\, \dd s+2\grr{\hat L_\vr} L_b\io\left|\kappa(x,z)-\kappa(y,z)\right|\, \dd z\,,
\end{align}
where $L_b$ is a constant depending on the  Lipschitz constants of $G$ and $G'$, $\|\kappa\|_{L^\infty(\Omega\times\Omega)}$,
$|\Omega|$, and $T$.
From \eqref{grad1}, using the assumptions $\chi_0\in\mathbf{V}$ and $\kappa\in W^{1,\infty}(\Omega\times\Omega)$ (cf.
Hypothesis~\ref{hyp1}(iii)), we immediately deduce that
\begin{equation}\no
|\nabla\chi(\cdot, t)|\leq C_\vr\left(1+\itt\left(|\nabla\teta(\cdot, s)|+|\nabla\chi(\cdot, s)|\right)\, \dd s\right)\quad \hbox{a.e.~in }\Omega\,.
\end{equation}
Now, with the help of Gronwall's lemma, we infer that
\begin{equation}\label{grad2}
|\nabla\chi(\cdot, t)|\leq C_\vr\left(1+\itt|\nabla\teta(\cdot, s)|\, \dd s\right)\quad \hbox{a.e.~in }\Omega\,.
\end{equation}
Using finally \eqref{grad2} with \eqref{e1fin1}, we get the desired estimate
\begin{equation}\label{egradchi}
\|\chi_n\|_{L^\infty(0,T;\mathbf{V})}\leq C_\vr\,,
\end{equation}
where $C_\vr$ denotes a positive constant depending increasingly on $\vr$.
From the definition \blu{\eqref{defutau}} of $u$, it also follows that
\begin{equation}\label{egradu}
\|u_n\|_{L^2(0,T; V)}\leq C_\vr\,.
\end{equation}


\subsection{Lower and upper bounds on $\teta$}\label{lowup}

In this \grr{subsection, we} first prove a bound for $\log\teta$ entailing the strict positivity of the absolute temperature (in the
limit when $n\to\infty$). Then,  we prove a (time dependent)  upper bound holding true for the solution component
$\teta_n$ for $n$ fixed, which
enables us to
proceed with the Moser iteration procedure in order to prove a uniform (independent of time, \grr{of $n$}, and of $\vr$) upper bound
on $\teta$. This permits us to remove the truncation parameter and to conclude the
existence proof. Finally, we will prove a lower bound (independent of $n$) on $\teta$ holding true under the additional Hypothesis~\ref{hypuni}
that we will use for the proof of uniqueness of solutions.

\paragraph{\bf \blu{Estimate on $\log\teta$.}}\label{tetapos}

Let us rewrite equation \eqref{p1tau}, by using \eqref{p2tau}, in the following form, for all
$z\in V$,
\begin{align}\label{p12}
&\io\partial_t\left(\blu{\frac1n\teta}+e(\teta,\chi)\right)\, z\dd x
+\io \grr{k(\bar\teta_n, \bar\chi_n)}\nabla\teta\cdot\nabla z\dd x
+\int_{\partial\Omega}\gamma(\teta-\blu{\teta_{\Gamma,n}})\,z\dd A\\
\no
&=\io \widetilde{\mu}_\vr(\teta)\chi_t^2+\teta \chi_t R(\teta,\chi)z\, \dd x\,,
\end{align}
where
$$R(\teta,\chi):=\sigma'(\chi)-\grr{s^\vr_{\chi}}(\teta,\chi)+\xi, \quad \xi\in \dvp(\chi), \quad |\xi(x,t)|\leq C\,\,\hbox{a.e.}\,,$$
$C$ being defined in Hypothesis~\ref{hyp1}(xi).
We prove now an estimate on $\log\teta$ in $L^2(0,T;V)$ by taking in \eqref{p1tau} $\displaystyle z=T(\teta)$,
where
\begin{equation}\label{defTdelta}
T(\teta):=-\left(1-\grr{\frac{1}{\teta}}\right)^-=
\begin{cases}
\displaystyle 1-\frac{1}{\teta}&\hbox{ for }\teta\leq 1\\
0&\hbox{ for }\teta\geq 1\,
\end{cases}\,.
\end{equation}
Notice that we are allowed to perform this estimate, \blu{with fixed $n$},
because \grr{$\teta\geq\e_n>0$ a.e.} for all $n$ (cf. Lemma~\ref{exiptau}).
We get, using equation \eqref{p2tau} and Hypothesis~\ref{hyp1}(iv),
\begin{align}\label{numera0}
&\frac{\dd}{\dd t} E(\teta,\chi)-\io \chi_t\int_{1}^\teta (c_V)_\chi(\xi,\chi) T(\xi)\, \dd\xi\dd x
+k_0\|\nabla(\log\teta)^-\|_H^2
\\
\no
&+\int_{\partial\Omega}\left(\teta-\grr{\teta_{\Gamma,n}}\right) T(\teta)\dd A\leq\io\widetilde{\mu}_\vr(\teta)|\chi_t|^2T(\teta)\dd x
+\io\teta T(\teta)\chi_t R(\teta,\chi)\dd x\,,
\end{align}
where
$$\displaystyle E(\teta,\chi)=-\io\int_{1}^\teta c_V(\xi,\chi)\left(1-\frac{1}{\xi}\right)^-\, \dd\xi\geq 0\,.$$
Note that the first term on the right-hand side of \eqref{numera0}
is nonpositive, while the other term can be estimated using estimates \eqref{estisro}, \eqref{estichi}, and \eqref{e1fin1} on
our solution $(\teta,\chi)$. Regarding the second term on the left-hand side in \eqref{numera0},
using \eqref{c4} in Hypothesis~\ref{hyp1}(vi), we obtain that
$$\left|-\io\chi_t\int_{1}^\teta (c_V)_\chi(\xi,\chi) T(\xi)\, \dd\xi\dd x\right|
\leq c_1 E(\teta,\chi)\|\dt \chi\|_{L^\infty(\Omega)^d}.$$
Moreover, we treat the boundary integral in the following way:
\begin{equation}\label{psi1}
\int_{\partial\Omega} (\teta-\grr{\teta_{\Gamma,n}})T(\teta)\dd A\geq \int_{\partial\Omega}\Psi(\teta)-\Psi(\grr{\teta_{\Gamma,n}})\dd A,
\end{equation}
where $\Psi$ is defined as
$$
\Psi(\teta):=
\begin{cases}
\displaystyle \teta-\log\teta&\hbox{ for }\teta\leq 1\\
1&\hbox{ for }\teta\geq 1\,
\end{cases}\,,
$$
and is a convex function on $[0,+\infty)$ such that $\Psi'(\teta)=T(\teta)$ for all $\teta\in [0,+\infty)$.
Using estimates \eqref{estichi} and \eqref{e1fin1}, and Hypothesis~\ref{hyp1}(x), we deduce that
\begin{align}\label{estpos0}
&E(\teta(x,t),\chi(x,t))+k_0\itt \|\grr{\nabla(\log\teta)^-}\|_H^2\dd \xi+\itt\int_{\partial\Omega}\Psi(\teta)\dd A\dd \xi\\
\no
&
\leq C_\vr\left(1+\itt E(\teta,\chi)\right)+E(\teta_{0,n},\chi_0)+\itt\int_{\partial\Omega}\Psi(\grr{\teta_{\Gamma,n}})\dd A\dd \xi
\end{align}
for a.e. $(x,t)\in Q_T$.
Now, in view of \eqref{s1}, we have that
\begin{align}\label{estiteta0}
E(\teta_{0,n},\chi_0)&=-\io\int_{1}^{\teta_{0,n}}c_V(\xi,\chi_0)\left(1-\frac{1}{\xi}\right)^-\, \dd\xi\dd x\\
\no
&\leq \io\int_0^{1}
c_V(\xi,\chi_0)\left(\frac{1}{\xi}-1\right)\, \dd\xi\dd x\leq \io s(1,\chi_0)\leq c_1|\Omega|\,.
\end{align}
Moreover, we can estimate the term $\itt\int_{\partial\Omega}\Psi(\grr{\teta_{\Gamma,n}})\dd A\dd \tau$ using Hypothesis~\ref{hyp1}(x) as follows:
\begin{equation}\label{estitetagamma}
\itt\int_{\partial\Omega}\Psi(\grr{\teta_{\Gamma,n}})\dd A\dd \xi\leq C\left(\|\grr{\teta_{\Gamma,n}}\|_{L^\infty(\Sigma_t)}
+\|\left(\log(\teta_{\Gamma, n})\right)^-\|_{L^1(\Sigma_t)}\right)\leq C\,.
\end{equation}
Using a standard Gronwall's lemma in \eqref{estpos0}, together with the estimates \eqref{estitetagamma} and \eqref{estiteta0},
we obtain
the desired bound
\begin{equation}\no
\|(\log\teta_n)^-\|_{L^2(0,T;V)}\leq C_\vr,
\end{equation}
which, together with estimate \eqref{e1fin1} for $\teta$ in $\LDV$, gives
\begin{equation}\label{estilog}
\|\log\teta_n\|_{L^2(0,T;V)}\leq C_\vr.
\end{equation}

\paragraph{\bf \blu{Lower bound for $\teta$ under Hypothesis~\ref{hypuni} (iii)--(vi)}.}\label{lowerbouuni}

\grr{Let us consider relation \eqref{p12}.} Notice that, \blu{by the previous estimates}, it follows that there exists a
positive constant \grr{$R_\vr$} such that $|R(\teta,\chi)|\leq \grr{R_\vr}$ a.e. in \grr{$Q_T$}.
Hence, \grr{for every $z \in V$ such that $z\geq 0$ a.e., we obtain the following inequality}
(notice that \grr{by virtue of \blu{Lemma~\ref{exiptau}}, we have here} $\teta_t\in L^2(0,T;H)$):
\begin{align}\label{ineqlow1}
&\io \left(\blu{\frac1n\teta_t}+c_V(\teta,\chi)\teta_t\right)\, z\dd x+
\io \grr{k(\bar\teta_n, \bar\chi_n)}\nabla\teta\cdot\nabla z\dd x
\\ \no
&=\io \grr{\tilde\mu_\vr(\teta)}\left(\chi_t+\frac{R(\teta,\chi)\teta}{2\grr{\tilde\mu_\vr(\teta)}}\right)^2
-\frac{R^2(\teta,\chi)\teta^2}{4\grr{\tilde\mu_\vr(\teta)}}z\, \dd x
\geq -\io\frac{\grr{R_\vr^2}\teta^2}{4\grr{\tilde\mu_\vr(\teta)}}z\, 
\dd x\,.
\end{align}
We now compare this inequality with the following ODE:
\begin{align}\label{eqw}
\tilde c(w) w_t=-\frac{\grr{R_\vr^2} w^2}{4\grr{\tilde\mu_\vr(w)}}, \quad w(0)=w_0\,,
\end{align}
where $\tilde c$ is defined in \grr{Hypothesis~\eqref{hypuni}(iii)}
and $w_0=\min_{x\in \Omega}\teta_0(x)\geq \teta_*$ 
(cf. Hypothesis~\ref{hypuni}(v)).

Notice that the solution $w$ is decreasing and does not vanish in finite time, due
to Hypothesis~\ref{hypuni}(iii).
\grr{The function $w$ does not depend on $x$, hence we may}
add to the ODE in \eqref{eqw} the term \grr{$-\dive(k(\bar\teta_n, \bar\chi_n)\nabla w)$, which is equal to $0$.
Using the fact that $w_t < 0$ and $\tilde c(w)\leq \frac{1}{n}+ c_V(w,\chi)$, we obtain,
subtracting \eqref{ineqlow1} from \eqref{eqw}, the inequality}
\begin{align}\no
&\io \left(\left(\blu{\frac1n}+c_V (w, \chi)\right) w_t
-\left(\blu{\frac1n}+c_V(\teta,\chi)\right)\teta_t\right)\, z\dd x
\\
\no
&{} +\io  k(\bar\teta_n, \bar\chi_n)\nabla(w-\teta)\cdot\nabla z\dd x
\ \leq\ \grr{\frac{R_\vr^2}{4}\io\left(\frac{\teta^2}
{\grr{\tilde\mu_\vr(\teta)}}-\frac{w^2}{\grr{\tilde\mu_\vr(w)}}\right)z\,
\dd x}\,.
\end{align}
\grr{We now take} as test function $z=H_\e(w-\teta)$, where $H_\e$
is the regularization of the Heaviside function $H$,
\begin{equation}\label{heaviside}
H_\e(v)=\begin{cases}
0 &\quad\hbox{if } v\leq 0\\
v/\e & \quad \hbox{if } v\in (0,\e)\\
1 & \quad\hbox{if } v\geq \e
\end{cases}\,.
\end{equation}
\grr{By virtue of  Hypothesis~\ref{hypuni} (iv)--(v)}, we deduce
\begin{equation}\no
\io \left(\left(\blu{\frac1n}+c_V(w, \chi)\right)w_t-\left(\blu{\frac1n}+c_V(\teta, \chi)\right) \teta_t\right) H_\e(w-\teta)\dd x\leq 0,
\end{equation}
and we can pass to the limit in this inequality for $\e\searrow 0$, getting
\begin{equation}\no
\io \left(\left(\blu{\frac1n}+c_V(w, \chi)\right)w_t-\left(\blu{\frac1n}+c_V(\teta, \chi)\right) \teta_t\right) H(w-\teta)\dd x\leq 0\,,
\end{equation}
that is,
\begin{equation}\label{ineqlow2}
\frac{\partial}{\partial t} \io \left(\left(\blu{\frac1n w}+e(w,\chi)\right)
-\left(\blu{\frac1n\teta}+e(\teta,\chi)\right)\right)^+\dd x\leq \io
\left(e_\chi(w,\chi)-e_\chi(\teta,\chi)\right)\chi_tH(w-\teta)\dd x\,.
\end{equation}
Notice now that, by  Hypothesis~\ref{hyp1}(vi) (cf. \eqref{c1} and \eqref{c4}), we have
$$ \left|e_\chi(w,\chi)-e_\chi(\teta,\chi)\right|\leq \max_{\tau\leq w_0} (c_V(\tau,\chi))|w-\teta|\leq \bar c|w-\teta|\,.$$
Integrating \eqref{ineqlow2} over $(0,t)$, and using the boundedness of $\chi_t$ in $L^\infty(Q_T)$,
we obtain, \grr{by
the choice of the initial data $\teta_0$ and $w_0$, that}
\begin{equation}\label{ee1}
\io \left(\left(\blu{\frac1n w} +e(w,\chi)\right)-\left(\blu{\frac1n\teta}+e(\teta,\chi)\right)\right)^+(t)\dd x\\
\ \leq \ \bar c\io (w-\teta)^+\dd x\,.
\end{equation}
\grr{For $w>\teta$, we have
\begin{align}\no
\frac{e(w,\chi)-e(\teta,\chi)}{w-\teta}=\frac{\int_\teta^w c_V (\tau, \chi)\dd \tau}{w-\teta}
\geq\frac{\int_\teta^w \tilde c (\tau)\dd \tau}{w-\teta}\geq \frac{1}{w}\int_0^w \tilde c(\tau)\dd \tau =:\tilde C(w)\,.
\end{align}
The function $\tilde C$ is nondecreasing, $\tilde C(w) > 0$ for $w>0$. Hence,
\begin{equation}\label{ineqe}
\left(\blu{\frac1n w}+e(w,\chi)\right)-\left(\blu{\frac1n\teta}+e(\teta,\chi)\right)
\geq \tilde C(w) (w-\teta)\quad \hbox{for } w\geq \teta\,.
\end{equation}
\grr{Inequality \eqref{ee1} then yields}
$$\grr{\tilde C(w(t))} \io (w-\teta)^+(t)\dd x\leq C\itt\io (w-\teta)^+\dd x \dd s$$
for every $t \in (0,T)$.} From Gronwall's lemma we conclude that
\begin{equation}\label{eslowbou}
\grr{\teta_n(x,t)\geq w(t)\quad\hbox{ a.e. in } Q_T\,.}
\end{equation}

\paragraph{\bf \blu{Upper bound for $\teta$.}}\label{upperboupar}

Let us denote the right-hand side in \eqref{p1tau} by
$$M(\teta,\chi):=\left((\lambda'(\chi)(t)+b[\chi](t))\dt\chi(t)
+\beta\dt\left(\vp(\chi(t))\right)\right)$$
which, due to the previous estimates is bounded by a positive constant, say, $\grr{\tilde M_\vr}$.
Then, we \grr{compare} the inequality, for all $z\in V$, $z\geq 0$ \grr{a.e.},
\begin{align}\label{p1tauineq}
&\io\dt\left({\frac{1}{n}\teta(t)}+{e}(\teta(t),{\chi(t)})\right)\, z\, \dd x
+\io {k(\bar\teta_n(t),\bar \chi_n(t))}\nabla\teta(t)\cdot\nabla z\dd x
\\
\no
& +\int_{\partial\Omega}\gamma(\teta(t)- {\teta_{\Gamma, n}(t)})z\dd A \leq \grr{\io\tilde M_\vr}\, z\dd x,
\end{align}
with the following ODE:
\begin{align}\label{eqwupper}
\grr{\frac1n}\, \dot{v}_n=\tilde M_\vr, \quad v(0)=v_0\,,
\end{align}
where \grr{$v_0 = \max\{\sup \teta_0, \sup\theta_\Gamma\}$ (cf. Hypothesis~\ref{hyp1} (viii) (x))}.

Then we have that
$$\grr{v_n}(t)=v_0+\grr{{\tilde M}_\vr n} t\,.$$
\grr{We proceed as above, adding to the ODE the term $-\dive (k(\bar\teta_n,\bar\chi_n)\nabla v_n)$, which is equal to $0$, and subtracting \eqref{eqwupper} from
\eqref{p1tauineq}. We} test the resulting inequality by $H_\e(\teta-v_n)$ (cf.~\eqref{heaviside}).
Using the fact that $\left(1/n+c_V(w,\chi)\right)\dot v_n\geq 0$ and that $\teta\in H^1(0,T;H)$, we can let $\e$ tend to $0$, getting
\begin{align}\no
&\io\frac{\partial}{\partial t}\left(\left(\frac1n\teta+e(\teta,\chi)\right)-\left(\frac1n v+e(v,\chi)\right)\right)^+\dd x\\
\no
&-\io\left(e_\chi(\teta,\chi)-e_\chi(v,\chi)\right)\chi_t H(\teta-v)\dd x\leq 0\,.
\end{align}
Using now the Lipschitz continuity of $e_\chi$ (cf. Hypothesis~\ref{hyp1}(vi)) \grr{and the boundedness of $\chi_t$}, we obtain that
\begin{align}\no
&\io\frac{\partial}{\partial t}\left(\left(\frac1n\teta+e(\teta,\chi)\right)-\left(\frac1n v+e(v,\chi)\right)\right)^+\dd x
\leq \grr{C_\vr}\io(\teta-v)^+\dd x\,.
\end{align}
Integrating over $(0,t)$ and using the choice of the initial condition \grr{$v_0$}, we infer
$$\io(\teta-v)^+(t)\dd x\leq \grr{C_\vr}\itt\io (\teta-v)^+\dd x\dd \tau,$$
and, applying Gronwall's lemma, we get the desired upper bound
\begin{equation}\label{upperbou}
\grr{\teta_n(x,t)\leq v_n(t)}\quad\hbox{a.e.~in }Q_T\,.
\end{equation}

\paragraph{\blu{Moser estimate.}}\label{moser}

In order to conclude the proof of existence of solutions to (\ref{p1})--(\ref{p3}), it remains only to prove that
the $\teta$ component of the solution $(\teta,\chi)$  is bounded from above independently
of $T$, $\vr$ and $n$.

This will enable us to choose $\vr$ sufficiently large in such a way that in this range
of values of $\teta$ we have $s_\chi=s_\chi^\vr$. To this end, we perform the following Moser estimate.

We will make repeated
use of the  well-known interpolation inequality \red{(cf. \cite{bin1})}
\begin{equation}\label{gn}
\|v\|_H\leq A\left(\eta\|\nabla v\|_H+\eta^{-N/2}\|v\|_{L^1(\Omega)}\right),
\end{equation}
which holds for every $v\in V$ and every $\eta\in (0,1)$,
with a positive constant $A$ independent of $v$ and $\eta$.

Following the ideas already exploited in \cite[Prop.~3.10, p.~296]{krs2},
for $j\in \nat$, we choose in \eqref{p1tau}
$z=((\teta_n-\teta_R)^+)^{\tk-1}\in L^2(0,T;H^1(\Omega))$, with $\teta_R=\max\{\Theta_\Gamma-1, \Theta, 1\}$, where
$$|\grr{\teta_{\Gamma, n}}(x,t)|\leq\Theta_\Gamma\quad\hbox{a.e. in }\grr{\Sigma_T}, \quad |\grr{\theta_{0,n}}(x)|\leq \Theta\quad\hbox{a.e. in }\Omega.$$
We know that $z\in L^2(0,T;H^1(\Omega))$, due to the upper bound on $\teta_n$ proved in \eqref{upperbou}.

Here below, we denote by $C_i,$ $i=1, 2,\dots$ some positive constants that may depend on the data of the
problem, but not on $j$, $T$, $\vr$, and $n$. We omit again the indices $n$ in $\teta_n$, $\chi_n$, for simplicity.

Now set $u=(\theta-\theta_R)$ and
take $z=(u^+)^{\tk-1}$ in  \eqref{p1tau} to  obtain that
\begin{align}\label{moser0}
& \left\langle\frac{1}{n}\teta_t(t)+(e(\teta(t), \chi(t)))_t,(u^+)^{\tk-1}\right\rangle
+\io k(\bar\teta_n(t),\bar\chi_n(t))\nabla\teta(t)\cdot\nabla (u^+)^{\tk-1}\dd x\\
\no
&\quad+\int_{\partial\Omega}\gamma(\teta(t)-\grr{\teta_{\Gamma,n}})\,(u^+)^{\tk-1}\dd A\\
\no
&=
-\io\left(\lambda'(\chi)(t)\chi_t(t)
+\beta\left(\vp(\chi(t))_t\right)-b[\chi](t)\chi_t(t)\right)\,(u^+)^{\tk-1}\dd x.
\end{align}
Our aim is to prove that there exists a positive constant $C^*$ (independent of $T$, $\vr$, and $n$) such that
\begin{equation}\label{moseresti}
\|\teta(t)\|_{L^\infty(\Omega)}\leq C^*\left(1+\log\vr\right)^{4+2N}\quad \hbox{for a.e.~}t\in (0, T).
\end{equation}
The first term on the left-hand side can be rewritten as
\begin{align}\label{moser1}
&\left\langle\frac{1}{n}\teta_t(t)+(e(\teta(t), \chi(t)))_t, (u^+)^{\tk-1}\right\rangle=
\frac1n\duav{\teta_t,(u^+)^{\tk-1}}+ \duav{\teta_t, c_V(\teta,\chi) (u^+)^{\tk-1}}\\
\no
&+\duav{e_\chi(\teta,\chi)\chi_t,(u^+)^{\tk-1}}.
\end{align}

Let us deal with the \blu{second} term on the right-hand side in \eqref{moser1},
using Hypothesis~\ref{hyp1}(vi) (notice that by \eqref{c2}, we have
 $c_V \geq \underline c$ in the set where $\theta\geq \teta_R\geq 1$) as follows:
\begin{equation}\label{moser2}
\duav{\teta_t, c_V(\teta,\chi) (u^+)^{\tk-1}}=\frac{\dd}{\dd t}E_j(u,\chi)
-\io\left(\int_{0}^u (c_V)_\chi(\xi+\teta_R,\chi)(\xi^+)^{\tk-1}\, \dd\xi\right)\chi_t\dd x,
\end{equation}
where
\begin{equation}\label{bouEP}
 \underline c \,2^{-j}\io (u^+)^{\tk}\dd x\leq E_j(u,\chi):=\io\int_{0}^u c_V(\xi+\teta_R,\chi)(\xi^+)^{\tk-1}\, \dd\xi\dd x
 \leq \bar{c}\, 2^{-j}\io (u^+)^{\tk}\dd x.
\end{equation}
Then, using \eqref{estichi} and Hypothesis~\ref{hyp1}(vi), we infer that
\begin{equation}
\io\left(\int_{0}^u (c_V)_\chi(\xi+\teta_R,\chi)(\xi^+)^{\tk-1}\, \dd\xi\right)\chi_t\dd x\leq C_1 (1+\log\vr)^2 E_j(u,\chi).
\end{equation}
Moreover,  by Hypothesis~\ref{hyp1}, we obtain the inequality
\begin{align}\nonumber
&\hspace{-15mm}\io k(\bar\teta_n,\bar\chi_n)\nabla u \nabla\left((u^+)^{\tk-1}\right)\,\dd x
+\int_{\partial\Omega}\gamma\,(u-\theta_{\Gamma,n}+\theta_R)\,(u^+)^{\tk-1}\,\dd A
\\ \nonumber
&\geq  k_0\frac{\tk-1}{2^{2j-2}}\io
\left|\nabla\left((u^+)^{\tko}\right)\right|^2\,\dd x
+\int_{\partial\Omega} \gamma\,\left((u^+)^{\tk}-(u^+)^{\tk-1}\right)\,\dd A\,.
\end{align}
Now, set $\pk=(u^+)^{\tko}$. Regarding the terms on the right-hand side in \eqref{moser0}
and the last term in \eqref{moser1}, using \eqref{estichi} and Hypothesis~\ref{hyp1}(vi), we realize that
\begin{align}\no
&-\io\left(\lambda'(\chi)(t)\chi_t(t)
+\beta\left(\vp(\chi(t))\right)_t-b[\chi](t)\chi_t(t)-e_\chi(\teta,\chi)\chi_t\right)\,(u^+)^{\tk-1}\dd x\\
\no
&\leq C_2 (1+\log\vr)^2\left(\io|\pk|^2\,\dd x+1\right)\,.
\end{align}
Let us now set $\displaystyle E_j^n=E_j+\frac{2^{-j}}{n}\io |\Phi_j|^2\dd x$.
Then, with the help of H\"older's
and Young's inequalities, we deduce that
\begin{align}\no
&\hspace{-10mm}
\frac{\dd}{\dd t}E_j^n(u,\chi)+\frac{k_0(\tk-1)}{2^{2j-2}}
\io|\nabla\pk|^2\,\dd x +\int_{\partial\Omega}\gamma\,|\pk|^2\,\dd A\\
\no
&\leq \left(1-\tkm\right)\int_{\partial\Omega}\gamma\,|\pk|^2\,\dd A
+\tkm\int_{\partial\Omega}\gamma\,\dd A
+C_1(1+\log\vr)^2E_j(u,\chi)\\
\no
&\qquad +C_2(1+\log\vr)^2\io\left(|\pk|^2+1\right)\,\dd x\,.
\end{align}
Multiplying the above inequality by $\tk$, in view of the upper bound for $E_j(u,\chi)$ in \eqref{bouEP}, we find out that
\begin{align}\label{phik}
&\hspace{-10mm}2^j\frac{\dd}{\dd t}E_j^n(u,\chi) + 2k_0\io|\nabla\pk|^2\,\dd x+\int_{\partial\Omega}\gamma|\pk|^2\, \dd A\\
\no
&\leq\,
2^{j}C_3(1+\log\vr)^2\left(1+\io |\pk|^2\,\dd x\right).
\end{align}
We now use the interpolation inequality (\ref{gn}) and note that, thanks to estimate \eqref{e1fin1}, we have
$$
\|\Phi_1\|_{L^1(\Omega)}^2=\left(\io u^+\dd x\right)^2\leq C_4(1+\log\vr)^4,\quad \|\Phi_j\|_{L^1(\Omega)}^2=\|\Phi_{j-1}\|_{H}^4\,.
$$
Thus, we derive the inequalities
\begin{align}\label{inter1}
\io|\Phi_1|^2\,\dd x &\leq\ 2A^2\left(\eta^2\|\nabla\Phi_1\|_H^2+\eta^{-N}
C_{4}(1+\log\vr)^4\right)\,,\\[2mm]\label{inter}
\io|\pk|^2\,\dd x &\leq\ 2A^2\left(\eta^2\|\nabla\Phi_j\|_H^2+\eta^{-N}\|
\Phi_{j-1}\|^4_H\right)\ \ \for j>1\,.
\end{align}
For $j=1$, we infer from \eqref{phik} and \eqref{inter1} that
\begin{align}\no
&\hspace{-5mm} 2\frac{\dd}{\dd t}E_1^n(u,\chi)
+2k_0\io|\nabla\Phi_1|^2\,\dd x+\int_{\partial\Omega}\gamma|\Phi_1|^2\, \dd A\\
\no
&\leq\ C_{5}(1+\log\vr)^2\,
\left(1+\eta^2\io|\nabla\Phi_1|^2\,\dd x + \eta^{-N}(1+\log\vr)^4\right).
\end{align}
Choosing $\eta=\sqrt{k_0}/(\sqrt{C_{5}(1+\log\vr)^2})$, we find that
\begin{equation}\label{phi1}
2\frac{\dd}{\dd t}E_1^n(u,\chi)+k_0\io|\nabla\Phi_1|^2\,\dd x+\int_{\partial\Omega}\gamma|\Phi_1|^2\, \dd A
\leq C_{6}\,(1+\log\vr)^{6+N}\,.
\end{equation}
For $j>1$, we get
\begin{align}\no
&\hspace{-5mm} 2^j\frac{\dd}{\dd t}E_j^n(u,\chi)
+2k_0\io|\nabla\Phi_j|^2\,\dd x+\int_{\partial\Omega}\gamma|\Phi_j|^2\, \dd A\\
\no
&\leq\ 2^jC_{7}(1+\log\vr)^2\,
\left(1+\eta^2\io|\nabla\Phi_j|^2\,\dd x + \eta^{-N} \|\Phi_{j-1}\|_H^4\right).
\end{align}
Choosing $\eta=\sqrt{k_0}/(\sqrt{2^{j}(C_{7}(1+\log\vr)^2)})$, we conclude
from \eqref{phik} and \eqref{inter} that
\begin{align}\label{phik2}
&\hspace{-5mm} 2^j\frac{\dd}{\dd t}E_j^n(u,\chi)+k_0\|\nabla \pk\|_H^2+\int_{\partial\Omega}\gamma|\pk|^2\,\dd A\\
\no
&\leq 2^{j\left(\frac{N}{2}+1\right)}C_{8}(1+\log\vr)^2\left(1+(1+\log\vr)^N\|\Phi_{j-1}\|_H^4\right)\,.
\end{align}
By assumption,  we have $\|\Phi_j(0)\|_H^2 = 0$. Hence, integrating \eqref{phi1} and \eqref{phik2} with respect to time
and using the lower bound in \eqref{bouEP}, we obtain that
\begin{align}\no
&\|\Phi_1(t)\|_H^2\leq C_{9}(1+\log\vr)^{6+N}\,,\\
\no
&\|\pk(t)\|_H^2\leq C_{10}2^{j\left(\frac{N}{2}+1\right)}(1+\log\vr)^2\left(1+(1+\log\vr)^N
\max_{0\leq\tau\leq t}\|\Phi_{j-1}(\tau)\|_H^4\right).
\end{align}
Define now
\[
z_j(t)=\max_{0\leq\tau\leq t}\sqrt{\|u(\tau)\|_{L^{2^j}(\Omega)}}
=\max_{0\leq\tau\leq t}\|\pk(\tau)\|_H^{2^{-j}}\,.
\]
Then we have
\begin{align}\no
&z_1(t)\leq C_{11}\,(1+\log\vr)^{(6+N)/4}\,,\\
\no
&z_j(t)\leq C_{10}^{2^{-j-1}} 2^{j(\frac{N}{2}+1)2^{-j-1}}(1+\log\vr)^{(N+2)2^{-j-1}}\max\{1,z_{j-1}(t)\}\,.
\end{align}
In particular, putting $y_j(t)=\max\{1,z_j(t)\}$, we get
\begin{align}\label{y1}
&y_1(t)\leq C_{11}(1+\log\vr)^{(6+N)/4},\\
\label{yj}
&y_j(t)\leq \left( C_{12}(1+\log\vr)^{(\frac{N}{2}+1)}\right)^{\tkm}
2^{\frac{j}{2}\,\left(\frac{N}{2}+1\right)\,\tkm}\, y_{j-1}(t)\,,\quad\hbox{for }j\geq 2.
\end{align}
Hence, passing to the logarithm in the inequality \eqref{yj} and summing up the result from 2 to $j$, we obtain
\begin{align}\no
\log y_j(t)&\leq\sum_{i=2}^j 2^{-i}\Big(
\log\left(C_{12}(1+\log\vr)^{\frac{N}{2}+1}\right) + \frac{i}{2}\,\left(\frac{N}{2}+1\right)\,\log 2\Big)+\log(y_1(t))\\
\no
&=\log\left(C_{12}(1+\log\vr)^{\frac{N}{2}+1}\right) \sum_{i=2}^j 2^{-i}+\left(\frac{N}{2}+1\right)\log 2\sum_{i=2}^j i 2^{-i}\\
\no
&
\quad+\log\left(C_{11}\,\left(1+\log\vr\right)^{\frac{(6+N)}{4}}\right)\\
\no
&\leq \log\left(2^{\left(\frac{N}{2}+1\right)C_{13}}\right)+\log\left(C_{14}\left(1+\log\vr\right)^{N+2}\right)\,,
\end{align}
independently of $j$ and $t>0$. Hence, we deduce
$$y_j(t)\leq 2^{\left(\frac{N}{2}+1\right)C_{13}}\left(C_{14}\left(1+\log\vr\right)^{N+2}\right), $$
independently of $j$ and $t>0$.
Choosing a proper
$\tilde{C}$, which is independent of $\vr$, we can conclude that
\begin{equation}\no
\sup_{t\geq 0,\,j\in\nat}\sqrt{\|u(t)\|_{L^{2^j}(\Omega)}}\ \leq\
\tilde{C}(1+\log\vr)^{2+N}\,.
\end{equation}
Formula (\ref{moseresti}) now immediately follows.


\subsection{Passage to the limit as $n\to\infty$}\label{limittau}

Our aim now is to pass to the
limit in (\ref{p1tau})--(\ref{initau}) as $n \to\infty.$
\smallskip

{}From \eqref{estichi}, \eqref{e1fin1}, \eqref{e1fin2}, \eqref{egradchi}, \eqref{egradu}, it follows
that, up to the extraction of some subsequence of $n$ as
$n\to\infty,$ there exist three functions $u,\, \teta: (0,T)\to
H$, $\chi: (0,T)\to \mathbf{H}$, such that we have
(as a consequence of the generalized
Ascoli theorem, see, e.g., \cite[Cor.~8, p.~90]{simon})
\begin{align}\label{cu}
u_n\to u\quad\,&\hbox{ weak star in }H^1(0,T;V')\cap
L^{\infty}(0,T;H)\cap L^2(0,T;V)\,\\
\no
&\hbox{ and strongly in }C^0([0,T];V')\cap L^2(0,T;H),\\
\label{cteta}
\teta_n\to \teta\quad\,&\hbox{ weak star in }L^\infty(0, T; H)\cap L^2(0,T;V),\\
\label{cchi}
\chi_n\to \chi\quad\,&\hbox{ weak star in }L^\infty(0, T; \mathbf{V}),\\
\label{ccchibis}
\dt\chi_n\to\dt\chi\quad\,&\hbox{ weak star in }L^\infty(\Omega\times (0,T))^d,\\
\label{cchiter}
\chi_n\to\chi\quad\,&\hbox{ strongly in }C^0([0,T];\mathbf{H})\,.
\end{align}
Moreover, as $\chi_n$ are uniformly bounded, it is easy to see that
\begin{align}
\label{cchitau}
\bar \chi_n\to \chi\quad&\hbox{ weak star in }L^\infty(0,T;\mathbf{V})\quad \hbox{and strongly in }L^q(Q_T),
\end{align}
for every $q\in (1,\infty)$, as $n\to\infty.$
Note that (cf. \eqref{defutau}) $u\geq 0$ and $\teta\geq 0$ a.e.
\blu{Then it turns out that,
at least for a subsequence, \grr{$u_n\to u$ a.e.} Now, we denote by $\psi(\,\cdot\,,\chi)$ the inverse function
of ${e}$ with respect to the first variable, that is, $e(\grr{\psi(w,\chi)})=w $ for all $(w,\chi)\in [0,\infty)\times {\cal D}(\vp)$.
Since $e$ is continuous, increasing in $\teta$, and such that $e(\teta,\chi)\geq e_1\teta-e_2$ for some constants $e_1, e_2>0$,
we infer that $\psi$ is continuous with a linear growth in $[0,\infty)\times {\cal D}(\vp)$ and $\psi(0,\chi)=0$. The Nemytskii operator is
therefore continuous in $L^2(Q_T)$,
and so this function is continuous and
increasing with a linear growth.} Hence, we have that
$$\red{\teta_n=\psi\left(\grr{u_n-\teta_n/n,\chi_n}\right)\to \psi(u,\chi)\quad\hbox{strongly in }L^2(0,T;H)\,.}$$
Hence, $\teta=\psi(u,\chi)$, or, equivalently,
$u=e(\teta,\chi)$. \blu{Finally, we check that also $\bar\teta_n$ converge strongly to $\teta$ in $L^2(0,T;H)$, at least
for a subsequence of $n\to\infty$.} Indeed, from the definition of $\bar\teta_n$, we get
\begin{align}\label{tetabar0}
&\int_0^T\io|\bar\teta_n(x,t)-\teta_n(x,t)|^2\, \dd x\dd t=\io\sum_{j=1}^n\int_{t_{j-1}}^{t_j}|\bar\teta_n(x,t)-\teta_n(x,t)|^2\dd t\dd x\\
\no
&=\int_0^{T/n}\|\teta_n(t)-\teta_{0,n}\|_H^2\dd t+\left(\frac{n}{T}\right)^2\sum_{j=2}^{n}\io\int_{t_{j-1}}^{t_j}\left|\int_{t_{j-2}}^{t_{j-1}}
\left(\teta_n(x,t)-\teta_n(x,\tau)\right)\dd\tau\right|^2\dd t\dd x\,.
\end{align}
By \eqref{cteta}, the first integral on the right-hand side of \eqref{tetabar0} can be estimated as follows:
\begin{align}\no
&\int_0^{T/n}\|\teta_n(t)-\teta_{0,n}\|_H^2\dd t\leq \frac{C}{n}\,.
\end{align}
Regarding the second term in \eqref{tetabar0}, we proceed  as follows:
\begin{align}\no
&\left(\frac{n}{T}\right)^2\sum_{j=2}^{n}\io\int_{t_{j-1}}^{t_j}\left|\int_{t_{j-2}}^{t_{j-1}}
\left(\teta_n(x,t)-\teta_n(x,\tau)\right)\dd\tau\right|^2\dd t\dd x\\
\no
&
\leq \frac{n}{T}\sum_{j=2}^{n}\int_{t_{j-1}}^{t_j}\int_{t_{j-2}}^{t_{j-1}}
\|\teta_n(t)-\teta_n(\tau)\|_H^2\dd\tau\dd t\\
\no
&\leq\frac{n}{T}\sum_{j=2}^{n}\int_{t_{j-1}}^{t_j}\int_{t-t_{j-1}}^{t-t_{j-2}}
\|\teta_n(t)-\teta_n(t-h)\|_H^2\dd h\dd t\\
\no
&\leq\frac{n}{T}\int_{0}^{2T/n}\left(\int_h^T
\|\teta_n(t)-\teta_n(t-h)\|_H^2\dd t\right)\dd h\,,
\end{align}
where we have used the new variable $h=t-\tau$. The last integral tends to $0$, because $\teta_n$ converge strongly to $\teta$,
which is mean continuous in $L^2(Q_T)$.
This implies that
$$
\bar\teta_n\to\teta\quad \hbox{strongly in }L^2(0,T;H),
$$
and this allows us to pass to the limit in (\ref{p1tau})--(\ref{initau}) as $n\to\infty$.
\blu{Notice moreover that from estimate \eqref{estilog} we also deduce that \grr{there exists a function $\zeta\,:\,(0,T)\to V$ such that}
$$
\log\teta_n\to \zeta\quad\hbox{weakly in }L^2(0,T;V).
$$
\grr{Using} the strong convergence of $\teta_n$ and the maximal monotonicity of the extended
$\log$ graph, we also deduce that $\zeta=\log\teta\in \LDV$, hence $\teta>0$ a.e.}

\subsection{Conclusion of the existence proof}\label{conclexi}

\blu{Let us now introduce the limit problem obtained by passing to the limit as $n\to\infty$ \grr{in (\ref{p1tau})--(\ref{initau}) in the
previous subsection.}}

\noindent {\sc Problem {(P)$_{\vr}$.}} For fixed $T>0$, and $\vr>0$,
find two functions $\teta\in  L^2(0,T;V)$ and
$\chi\in L^{\infty}(\Omega\times(0,T))^d$, $\chi_t\in L^{\infty}(\Omega\times(0,T))^d$
 such that,
for $t \in (0,T)$, we have
\begin{align}\label{p1piu}
& \duav{(e(\teta(t), \chi(t)))_t, z}
+\io k(\teta(t),\chi(t))\nabla\teta(t)\cdot\nabla z\dd x+\int_{\partial\Omega}\gamma(\teta(t)-\teta_\Gamma)\,z\dd A\\
\no
&=
-\io\left(\lambda'(\chi)(t)\chi_t(t)
+\beta\left(\vp(\chi(t))\right)_t-b[\chi](t)\chi_t(t)\right)\,z\dd x\quad\forall z\in V,\\
\label{p2piu}
&\grr{\tilde{\mu}_\vr}(\teta(t))\chi_t(t)
+(\beta+\teta(t))\dvp(\chi)(t)\ni
-\lambda'(\chi)(t)-\teta(t)\sigma'(\chi)(t)\\
\no
&-b[\chi](t)-e_\chi(\teta(t),\chi(t))
+\teta(t) s_\chi^\vr(\teta(t),\chi(t))\quad\hbox{a.e.~in }\Omega\,,
\end{align}
and for $t=0$ the functions $e(\teta,\chi)$ \grr{and} $\chi$ satisfy the initial conditions
\begin{equation}
\label{inipiu}
\displaystyle e(\teta,\chi)(0)=u_0,\quad \chi(0)=\chi_0\,.
\end{equation}

We thus have proved the following result.
\bepr\label{exirho} Let
Hypothesis~\ref{hyp1} hold true and let $\vr>0$, $\grr{\tilde{\mu}_\vr}$ and
$s_\chi^\vr$ be defined as at the beginning of
Section~\ref{approxi}. Then {\sc Problem {(P)$_{\vr}$}} has a
solution in $(0,T)$.
\empr

In order to conclude the proof of Theorem~\ref{thexi}, we only need to remove the truncation $\vr$. This can be done using
the Moser estimate \eqref{moseresti}.

This bound for $\teta$, indeed,  allows us, for a suitable $\vrh\geq 1$ satisfying
$$C^*(1+\log\vr)^{4+2N}\leq \vr/2,$$
with, e.g., $C^*=\tilde{C}^2$, to remove the truncation
$\vrh$ and to conclude that the solution to {\sc Problem (P$_\vr$)} is in fact a solution to {\sc Problem (P)}.

\grr{Finally, let us note that we have now found a solution on $(0,T)$ and we deduce that $\chi$ is weakly continuous
with values in $\mathbf{V}$, $e(\teta,\chi)\in C^0([0,T];H)$, and $\teta$ is bounded uniformly in time in $L^\infty(\Omega)$ due
to the Moser estimate. Hence, we can continue the solution
starting from time $T$ and extend it on the whole time interval $(0,\infty)$.}

This concludes the proof of Theorem~\ref{thexi}.\QED


\section{Uniqueness}\label{proofuni}

In this section, we prove Theorem~\ref{thuni}. Hence, we assume Hypothesis~\ref{hypuni}.
First let us note that the lower bound for $\teta$ \eqref{lowbou} directly follows by passing
to the limit as $n$ tends to $\infty$ in \eqref{eslowbou}.
Now we are ready to proceed with the proof of
uniqueness and finally we will prove the regularity result \eqref{regotetat} under the further assumption \eqref{regotetagamma}.

\paragraph{Uniqueness.} In what follows,
we denote by $R_0, R_1, R_2, \dots$ suitable constants that possibly depend on  $T$,
but not on the solutions. We start with the proof of uniqueness of solutions. Equation (\ref{p2})
is for (almost) all $x\in \Omega$ of the form
\eqref{diffinclu}, with
\begin{align}\no
\alpha(\theta) &= \frac{\mu(\theta)}{\beta + \theta}\,,\\ \no
g = \ell[\theta,\chi] &= -\frac{1}{\beta + \theta}\,
\left(
\lambda'(\chi)+\teta\sigma'(\chi)+b[\chi]+e_\chi(\teta,\chi)-\teta s_\chi(\teta,\chi)\right)\,.
\end{align}
Within the range $0<  \teta < \overline\teta$
and $\chi \in \dcvp$, $|\chi_t| \le C$, of admissible values for the solutions, and,
thanks to Hypothesis~\ref{hyp1} (ii) (iii) (v) (vi), all nonlinearities
in (\ref{p1})--(\ref{p2}) are Lipschitz continuous.
 Using the notation from Theorem \ref{thuni}, we obtain,
as a consequence of (\ref{cdg}), for a.\,e. $(x,t) \in Q_\infty$
the estimate
\red{
\begin{align}\label{uni1}
&\itt |\hat\chi_t(x,\tau)|\dd \tau+|\hat\chi(x,t)|\\
\no
&
\le R_0(T) \left(|\hat\chi_0(x)|+\itt\left(|\hat\theta(x,\tau)| + |\hat\chi(x,\tau)| + \io |\hat\chi(y,\tau)|\,\dd y
\right)\dd\tau\right),
\end{align}
with some constant $R_0(T)$}. Integrating over $\Omega$, and by Gronwall's argument, we obtain that
\begin{equation}\label{uni2}
\io |\hat\chi(y,t)|\,\dd y \le R_1
\left(\io |\hat\chi_0(y)|\,\dd y + \int_0^t\io |\hat\theta(y,\tau)|\,\dd y\,\dd\tau
\right)\,,
\end{equation}
and hence, we get
\begin{equation}\label{uni3}
\int_0^t|\hat\chi_t(x,\tau)|\,\dd\tau + |\hat\chi(x,t)|
\le R_2 \left(|\hat\chi_0(x)|+\int_0^t|\hat\theta(x,\tau)|\,\dd\tau
+  \int_0^t\io |\hat\theta(y,\tau)|\,\dd y\,\dd\tau\right)
\end{equation}
for a.\,e. $x\in\Omega$ and every $t\in [0,T]$. In particular,
\begin{equation}\label{uni5}
\int_0^t\io |\hat\chi_t(x,\tau)|\,\dd x\,\dd\tau \le R_3\left(
\io|\hat\chi_0(x)|\,\dd x
+ \int_0^t\io |\hat\theta(x,\tau)|\,\dd x\,\dd\tau\right)\,.
\end{equation}
We now multiply (\ref{uni3}) by $|\hat\chi(x,t)|$ and integrate over $\Omega$
to obtain for all $t\in [0,T]$ that
\begin{equation}\label{uni4}
\io |\hat\chi(x,t)|^2\,\dd x \le R_4\left( \|\hat\chi_0\|_{\mathbf{H}}^2 +
\int_0^t\io |\hat\theta(x,\tau)|^2\,\dd x\,\dd\tau\right)\,.
\end{equation}
The crucial point is to exploit Eq.~(\ref{p1}) properly.
Notice first that we have
\begin{equation}\label{uni6}
b[\chi]\chi_t(x,t) = 2B[\chi]_t(x,t) +
2\io \kappa(x,y)\,G'(\chi(x,t) - \chi(y,t))\,\chi_t(y,t)\,\dd y\,.
\end{equation}
We integrate the difference of the two equations (\ref{p1}),
written for $(\theta_1, \chi_1)$ and $(\theta_2, \chi_2)$, from $0$ to $t$,
rewriting the terms $b[\chi_i](\chi_i)_t$ according to (\ref{uni6}).
Take $z=K(\theta_1)-K(\teta_2)$ in the resulting equation, where $K(u)=\int_0^u \bar k(s)\,\dd s$, $u\in\RR$,
and integrate it again over $(0,t)$.
Using the lower bound for $\teta$, the Lipschitz continuity
of all nonlinearities ($\vp$ is Lipschitz continuous on $\dcvp$ with
constant $C$), the properties of $K$ (cf.~Hypothesis~\ref{hyp1}(iv)), 
and denoting
$$\hat\Theta(x,t) =
\int_0^t {\pier (K(\theta_1 - K(\theta_2))}(x,\tau) \dd\tau,
$$
using \eqref{eslowbou}, we obtain for each $t\in (0, T)$ that
\begin{align}\label{euni1}
&k_0c_V(w)\itt\io|\hat\theta(x,\tau)|^2\,\dd x\,\dd\tau +
\frac{1}{2} \io |\nabla\hat\Theta(x,t)|^2 \dd x
\\
\no
&\ \leq R_5\Big(
\|\hat\theta_0\|_H^2 + \|\hat\chi_0\|_{\mathbf{H}}^2 +
\itt\io |\hat\chi(x,\tau)|^2\,\dd x\,\dd\tau   \\
\no
&\hspace{0.5cm}+\itt\int_0^\xi \io\io
\kappa(x,y)|\hat\chi_t(y,\tau)| |\hat\theta(x,t)|\dd x\,\dd y\,\dd\tau\dd \xi\Big)\,.
\end{align}
The last term on the right-hand side of the above inequality
can be estimated, using~(\ref{uni5}),~by
\begin{align}\no
&\int_0^t \io\io \kappa(x,y)|\hat\chi_t(y,\tau)| |\hat\theta(x,t)|
\dd x\,\dd y\,\dd\tau\le R_6 \io|\hat\theta(x,t)|\dd x
\int_0^t \io |\hat\chi_t(y,\tau)|\dd y\,\dd\tau\\ \no
&\quad \le R_7 \left(\io|\hat\theta(x,t)|^2\dd x\right)^{1/2}
\left(\io|\hat\chi_0(x)|^2\,\dd x
+ \int_0^t\io |\hat\theta(x,\tau)|^2\,\dd x\,\dd\tau\right)^{1/2}\,.
\end{align}
Combining the last two inequalities again with the Gronwall's lemma, we 
obtain for each $t\in [0,T]$ the estimate
\begin{align}\label{uni7}
&\int_0^t \io |\hat\theta(x,\tau)|^2 \dd x\,\dd\tau
+ \io |\nabla\hat\Theta(x,t)|^2 \dd x 
\\ \no &\qquad 
\le R_8 \left(\|\hat\theta_0\|_H^2 + \|\hat\chi_0\|_{\mathbf{H}}^2
+ \int_0^t\io |\hat\chi(x,\tau)|^2\,\dd x\,\dd\tau \right)\,.
\end{align}
We now multiply (\ref{uni7}) by $2R_4$, add the result to (\ref{uni4}),
and see that Gronwall's argument can be applied again to arrive at
the final estimate
\begin{equation}\label{uni8}
\io |\hat\chi(x,t)|^2\dd x + \int_0^t \io |\hat\theta(x,\tau)|^2 \dd x\dd\tau
\le R_8 \left( {\pier \|\hat\theta_0\|_H^2 + 
\|\hat\chi_0\|_{\mathbf{H}}^2 }
\right).
\end{equation}

\paragraph{Regularity.} We prove now the regularity \eqref{regotetat} for $\teta$ under the further assumption
\eqref{regotetagamma}. In order to do that, let us consider, instead of \eqref{p1tau},
the following approximated equation
\begin{align}\label{p1nbis}
&\io\dt\left(\red{\frac{1}{n}\teta_n(t)}+{e}(\teta_n(t), \blu{\chi_n(t)})\right)\, z\, \dd x
+\io \red{\bar k(\teta_n(t))}\nabla\teta_n(t)\cdot\nabla z\dd x
\\
\no
&\qquad =
-\io\left((\lambda'(\chi_n)(t)+b[\chi_n](t))\dt\chi_n(t)
+\beta\dt\left(\vp(\chi_n(t))\right)\right)z\dd x\,.
\end{align}
Since (cf. Hypothesis~\ref{hypuni}(i))  the heat conductivity $\bar k$ is independent of $\chi$,
we can now test \eqref{p1nbis} by $K(\teta_n)_t$, where $K(\teta_n)=\int_0^{\teta_n} \bar k(s)\,\dd s$
and, integrating over $(0,t)$, we obtain, using \eqref{eslowbou}, the estimate
\begin{align}\no
& k_0 \left(\frac1n+c_V(w)\right)\itt \io |(\teta_n)_t |^2 \dd x\dd\tau
+\frac 12 \|\nabla(K(\teta_n))(t)\|_H^2
\\ \no
&\qquad \leq \frac12\|\nabla(K(\teta_{0n}))\|_H^2
+C k_1\itt\io|(\teta_n)_t|\dd x\dd \tau.
\end{align}
Here $C$ is a bound in $L^\infty(Q_T)$ for the term $(\lambda'(\chi_n)+b[\chi_n])\dt\chi_n
+\beta\dt\left(\vp(\chi_n)\right)$ that we already obtained in Theorem~\ref{thexi}, and $k_0$, $k_1$ are the positive constants
introduced in Hypothesis~\ref{hyp1}(iv). Using now assumption \eqref{regotetagamma} and {\pier the conclusions of Theorem~\ref{thexi},}
we obtain the (independent of $n$) bound
\begin{equation}\no
\|(\teta_n)_t\|_{\LDH}+\|K(\teta_n)\|_{\LIV}\leq C\,,
\end{equation}
which leads immediately (due to Hypothesis~\ref{hyp1}(iv)) to the desired estimate \eqref{regotetat}.

With this, Theorem \ref{thuni} is proved.
\QED

\beos\label{obc}
It seems to the authors that, in the case when $k$ depends only on $\theta$, 
if we replace the boundary condition \eqref{bou} with
\begin{equation}  
\label{bouK}
k(\teta)\nabla\teta\cdot{\bf n}+\gamma( K(\teta) -K(\teta_\Gamma)) =0\quad\hbox{on $\partial \Omega$},
\end{equation}
where $K$ is the primitive of $k$ as in this section, then the uniqueness
result holds for a general nonnegative $\gamma$. 
Moreover, the regularity estimates can be 
adapted in order to cover also this case:
the details of this proof will be developed in a subsequent paper. However,
let us point out that the regularity properties \eqref{regotetat}
remain true for more general $\gamma$ also in the case of 
boundary condition~\eqref{bou}:  the argument relies on a suitable use of the Gronwall lemma.
\eddos



\end{document}